\documentclass[12pt]{article}
\usepackage[T1]{fontenc}
\usepackage[utf8]{inputenc}
\usepackage{lmodern}
\usepackage[intlimits]{amsmath}
\usepackage{amssymb,amsthm,mathrsfs,mathtools,leftidx,bm}
\usepackage[margin=1in]{geometry}
\usepackage{indentfirst,enumitem,xcolor}
\usepackage[nottoc]{tocbibind}
\usepackage{microtype}
\usepackage[hidelinks]{hyperref}
\newcommand{\Mp}{\mathrm{Mp}}
\newcommand{\mpdr}{\Mp(n,\R)}
\newcommand{\RMp}{R^{\Mp}}
\newcommand{\ThetaMp}{\Theta^{\Mp}}
\newcommand{\KMp}{K_{\Mp}}

\theoremstyle{plain}                    
\newtheorem{theorem}{Theorem}[section]
\newtheorem*{theorem*}{Theorem}
\newtheorem{lemma}[theorem]{Lemma}    
\newtheorem{prop}[theorem]{Proposition}
\newtheorem{proposition}[theorem]{Proposition}

\theoremstyle{definition}

\theoremstyle{remark}
\newtheorem{remark}{Remark}[section]

\newcommand{\Z}{\mathbb{Z}}
\newcommand{\Q}{\mathbb{Q}}
\newcommand{\R}{\mathbb{R}}
\newcommand{\C}{\mathbb{C}}

\newcommand{\T}{\mathbb{T}}

\newcommand{\uhp}{\mathfrak H}

\newcommand{\ha}{\frac{1}{2}}
\newcommand{\tha}{\tfrac{1}{2}}

\newcommand{\matr}[4]{\left( \begin{matrix} #1 & #2 \\ #3 & #4 \end{matrix} \right) }

\newcommand{\smatr}[4]{\bigl( \begin{smallmatrix} #1 & #2 \\ #3 & #4 \end{smallmatrix} \bigr) }

\newcommand{\inv}{^{-1}}
\newcommand{\diag}{\mathrm{diag}}

\newcommand{\bs}{\backslash}

\newcommand{\lp}{\left(}
\newcommand{\rp}{\right)}
\newcommand{\labs}{\left|}
\newcommand{\rabs}{\right|}

\newcommand{\lcur}{\left\{}
\newcommand{\rcur}{\right\}}

\newcommand{\Heisd}{\mathbb{H}(\mathbb R^n)}

\newcommand{\LtRd}{L^2(\R^n)}
\newcommand{\Si}{\mathcal{S}}

\newcommand{\spdr}{\mathrm{Sp}(n, \R)}

\newcommand{\spdz}{\mathrm{Sp}(n, \Z)}

\newcommand{\sltr}{\mathrm{SL}(2, \R)}

\newcommand{\sltz}{\mathrm{SL}(2, \Z)}

\newcommand{\sltrn}{(\mathrm{SL}(2, \R))^n}
\newcommand{\sltzn}{(\mathrm{SL}(2, \Z))^n}
\newcommand{\AffSym}{\mathrm{ASp}(n,\R)}
\newcommand{\Un}{\mathrm{U}(n)}

\newcommand{\tg}{g}

\newcommand{\Jac}{\mathrm{J}_n(\R)}
\newcommand{\prodsltr}{(\sltr)^n}

\newcommand{\AffSL}{\R^{2n} \rtimes \prodsltr}

\newcommand{\frakn}{\mathfrak{n}}
\newcommand{\fraka}{\mathfrak a}
\newcommand{\calE}{\mathcal E}
\newcommand{\Gr}{\mathrm{Gr}}
\newcommand{\calG}{\mathcal{G}}
\newcommand{\oGam}{\overline{\Gamma}_\Theta}
\newcommand{\oDel}{\overline{\Delta}_\Theta}

\newcommand{\Sp}{\mathrm{Sp}}
\newcommand{\GL}{\mathrm{GL}}
\newcommand{\Diag}{\mathrm{Diag}}

\newcommand{\calF}{\mathcal F}

\newcommand{\calS}{\mathcal{S}}

\newcommand{\calD}{\mathcal{D}}

\newcommand{\vecxi}{\text{\boldmath$\xi$}}

\newcommand{\veca}{\text{\boldmath$a$}}
\newcommand{\vecb}{\text{\boldmath$b$}}
\newcommand{\vecc}{\text{\boldmath$c$}}
\newcommand{\vecp}{\text{\boldmath$p$}}

\newcommand{\vecv}{\text{\boldmath$v$}}
\newcommand{\vecu}{\text{\boldmath$u$}}
\newcommand{\vecw}{\text{\boldmath$w$}}
\newcommand{\vecm}{\text{\boldmath$m$}}
\newcommand{\vecn}{\text{\boldmath$n$}}
\newcommand{\vecx}{\text{\boldmath$x$}}
\newcommand{\vecy}{\text{\boldmath$y$}}
\newcommand{\vecz}{\text{\boldmath$z$}}
\newcommand{\vecr}{\text{\boldmath$r$}}
\newcommand{\vecs}{\text{\boldmath$s$}}

\newcommand{\vecnull}{\text{\boldmath$0$}}
\newcommand{\bn}{\text{\boldmath$0$}}
\newcommand{\vecphi}{\text{\boldmath$\phi$}}
\newcommand{\Sym}{\mathrm{Sym}}
\newcommand{\scrP}{\mathscr{P}}

\newcommand{\scrSP}{\mathscr{S\!P}}

\newcommand{\veps}{\varepsilon}

\newcommand{\Leb}{\mathrm{Leb}}
\newcommand{\vol}{\mathrm{Vol}}
\newcommand{\dlimmeas}{\nu_\Theta}

\newcommand{\de}{\:\mathrm d}
\newcommand{\ind}{\mathbf{1}}

\def\trans{\,^\mathrm{t}}
\def\d{\mathrm{d}}
\def\e{\mathrm{e}}
\def\i{\mathrm{i}}

\def\Re{\operatorname{Re}}
\def\Im{\operatorname{Im}}

\numberwithin{equation}{section}

\title{Limit laws for Siegel theta sums}
\author{Jory Griffin\thanks{{School of Mathematics, University of Bristol, Bristol BS8 1UG, U.K.}} \and  Jens Marklof\thanks{School of Mathematics, University of Bristol, Bristol BS8 1UG, U.K.; j.marklof@bristol.ac.uk} \and Tariq Osman\thanks{Institut für Mathematik, Universit\"at Zürich, CH-8057 Z\"urich, Switzerland; tariq.osman@math.uzh.ch}}

\date{24 November 2026}

\begin{document}
\maketitle

\begin{abstract}
\noindent
Siegel theta sums are exponential sums over the integer lattice whose phases are determined by a quadratic form. For a random choice of quadratic form, we prove a limit law for the modulus of the normalised sum as the summation range tends to infinity. The limiting distribution depends on the space from which the quadratic form is sampled. We prove precise heavy-tail asymptotics for random sampling from (a) the full space of quadratic forms and (b) the space of diagonal forms. This study extends classic results by Jurkat and Van Horne for one-variable theta sums and by G\"otze and Gordin for theta series on Siegel half-spaces. Our approach uses the Shale--Weil representation of the metaplectic group and the ergodic theory of unipotent flows.
\end{abstract}

\tableofcontents

\section{Introduction}

Let $\Sym_n(\R)$ denote the set of symmetric $n \times n$ matrices with real entries. The aim of this paper is to establish limit laws for sums of the form
\begin{equation}\label{eq:theta-sum}
S_f(M,X)=\sum_{\vecm\in\Z^n}f(M^{-1}\vecm)
e\!\lp\tha\vecm X\trans\vecm\rp,
\end{equation}
where $f:\R^n\to\C$ is a rapidly decreasing cut-off function, $M>0$, $e(x)=\e^{2\pi\i x}$, and $X\in\Sym_n(\R)$. We represent elements of $\R^n$ as row vectors. 
If $f(\vecv)=\e^{-\pi\vecv Y\trans\vecv}$ with $Y$ positive definite,
then
\begin{equation}\label{eq:gaussian-siegel-theta}
S_f(M,X)=\sum_{\vecm\in\Z^n}
e\!\lp\tha\vecm Z\trans\vecm\rp,
\qquad Z=X+\i M^{-2}Y.
\end{equation}
If instead $f=\ind_\Omega$ for a bounded set $\Omega\subset\R^n$, then \eqref{eq:theta-sum} is the finite exponential sum
\begin{equation}
  \label{def:finite-theta-sum}
  S_\Omega (M, X) := \sum_{\vecm \in \Z^n \cap M \Omega}e(\tha \vecm X \trans\vecm).
\end{equation}

One-variable theta sums have been studied since Hardy and Littlewood
\cite{Hardy-Littlewood1914partII}; see
\cite{Fiedler-Jurkat-Korner77,Flaminio-Forni2006,Fedotov-Klopp-2012,Cellarosi-Marklof,Cellarosi-Osman-rational-tails}
and references therein.  Much less is known in several variables.  Marklof and Welsh \cite{Marklof-Welsh-higher-rank-theta-sums} proved that, for $\varepsilon>0$, there exists a set $\mathcal{X}(\veps)$ of full Lebesgue measure such that for any $X \in \mathcal{X}(\veps)$ the following bound holds
\begin{equation}\label{upperJM}
S_f(M,X)=O_{f,X,\varepsilon}
\!\left(M^{n/2}(\log M)^{\frac1{2n+2}+\varepsilon}\right),
\end{equation}
for a broad class of cut-off functions $f$.
This improves an earlier bound of Cosentino and Flaminio
\cite{Cosentino-Flaminio}.  Bounds with an additional linear phase are
discussed in \cite{Fedotov-Klopp-2012,Marklof-Welsh-higher-rank-theta-sumsII, Cellarosi-Osman-special-bounds}.

We will mostly work with cut-off functions $f\in\Si(\R^n)$ of Schwartz class. Riemann-integrable cut-off functions $f$ are understood to be bounded, rapidly decreasing, and Riemann-integrable on every compact rectangle.

Throughout, we use bold lower-case letters to denote row vectors
and other ordered tuples.  If $\boldsymbol z=(z_1,\ldots,z_n)$, then
$\diag(\boldsymbol z)$ denotes the diagonal matrix whose $j$th diagonal
entry is $z_j$.  
Products of tuples with row vectors are understood componentwise and will be denoted by
\begin{equation}
\boldsymbol\varepsilon\vecx = (\varepsilon_1x_1,\ldots,\varepsilon_nx_n).
\end{equation}

Define the symmetry operators
\begin{equation}\label{eq:parity-projections}
(P_+f)(\vecx):=\frac{f(\vecx)+f(-\vecx)}2,
\qquad
(P_\square f)(\vecx):=2^{-n}
\sum_{\boldsymbol\varepsilon\in\{\pm1\}^n}
f(\boldsymbol\varepsilon\vecx).
\end{equation}

Then, pairing lattice points under the corresponding reflections gives
\begin{equation}\label{eq:parity-sums}
S_f(M,X)=S_{P_+f}(M,X),
\qquad
S_f(M,\diag(\boldsymbol x))
=S_{P_\square f}(M,\diag(\boldsymbol x)).
\end{equation}

We also note that the first sum in \eqref{eq:parity-sums} is periodic under translation by the lattice $\mathcal L_n :=\{L\in\Sym_n(\Z):L_{ii}\in2\Z\}$, and the diagonal sum is periodic under independent translations by $2$ in each coordinate, and so they are naturally defined over 
\begin{equation}\label{eq:sampling-tori}
\mathbb T_\Theta:=\mathcal{L}_n\bs\Sym_n(\R),\quad
\mathbb T_\Theta^\Diag:=\Diag_n(2\Z)\bs\Diag_n(\R).
\end{equation}
Let $m_\Theta$ and $m_\Theta^\Diag$ be the normalised Lebesgue measures on $\T_\Theta$ and $\T_\Theta^\Diag$, respectively. Thus $\d m_\Theta=2^{-n}\d X$ on the cuboid $\mathcal{C}_\Theta=[0,2)^n\times[0,1)^{n(n-1)/2}$ and $\d m_\Theta^\Diag=2^{-n}\d\vecx$ on $[0,2)^n$, where $\d X$ and $\d\vecx$ are the coordinate Lebesgue measures on these cuboids. 

Our first theorem generalises the limit laws for single-variable sums, first appearing in the works of Jurkat and Van Horne \cite{Jurkat-vanHorne1981-proofCLT,Jurkat-vanHorne1982-CLT,Jurkat-VanHorne1983} and Marklof \cite{Marklof-1999}. In particular, the asymptotic \eqref{tailasymptotic001} gives a nonzero power-law tail whenever $P_+f\neq0$. The limit theorem given in \eqref{LL001} was proven for $n \geq 2$ by G\"otze and Gordin in \cite{Gotze-Gordin-theta-limit-thm} for a Gaussian cut-off $f$, that is, for sums of the form \eqref{eq:gaussian-siegel-theta}, but there they do not address tail asymptotics or convergence of higher moments \cite[Remark 4.1]{Gotze-Gordin-theta-limit-thm}.

\begin{theorem}\label{thm:main1}
Let $n\geq2$ and $f\in\Si(\R^n)$. Let $\lambda \ll m_\Theta$ be a Borel probability measure on $\T_\Theta$. There is a right-continuous non-increasing function $\Psi_f:\R\to[0,1]$ with $\Psi_f(R)=1$ for $R<0$, which is independent of $\lambda$, with the following properties.
\begin{enumerate}
\renewcommand{\labelenumi}{\emph{(\roman{enumi})}}
\item \emph{Limit law.}
  At every continuity point of $\Psi_f$, we have
\begin{equation}\label{LL001}
\Psi_f(R)=\lim_{M\to\infty}\lambda\left\{ X\in\T_\Theta :M^{-n/2}|S_f(M,X)|>R\right\}.
\end{equation}

\item \emph{Tail asymptotic.}
For every $\varepsilon>0$,
\begin{equation}\label{tailasymptotic001}
\Psi_f(R)=C_{n,f}R^{-2(n+1)}
+O_{n,f,\varepsilon}\!\left(R^{-2(n+2)+\varepsilon}\right),
\qquad R\to\infty,
\end{equation}
where $C_{n,f}$ is given in \eqref{eq:generic-tail-constant}, and in particular  $C_{n,f}$ is positive if and only if $P_+f\neq0$.

\item \emph{Moments.}
If the density of $\lambda$ is in $L^\infty (\T_\Theta)$, then for  every $0<\beta<2(n+1)$,
\begin{equation}\label{eq:full-moment-convergence-intro}
\lim_{M\to\infty}\int_{\T_\Theta}
\left(M^{-n/2}|S_f(M,X)|\right)^\beta\,\d\lambda(X)
=\beta\int_0^\infty R^{\beta-1}\Psi_f(R)\,\d R<\infty.
\end{equation}
In particular, when $\beta=2$ the limit \eqref{eq:full-moment-convergence-intro} evaluates to  $2\|P_+f\|_2^2$.
\end{enumerate}
\end{theorem}

We contrast the general case with sampling only in the subspace of diagonal $n \times n$ matrices $\Diag_n(\R)$, which has zero Lebesgue measure in $\Sym_n(\R)$ when $n\geq2$.

\begin{theorem}\label{thm:main2}
Let $n\geq2$ and $f\in\Si(\R^n)$. Let $\lambda \ll m_\Theta^{\Diag}$ be a Borel probability measure on $\T_\Theta^{\Diag}$.  There is a right-continuous non-increasing function $\Psi_f^\Diag:\R\to[0,1]$ with $\Psi_f^\Diag(R)=1$ for $R<0$, which is independent of $\lambda$, with the following properties.
\begin{enumerate}
\renewcommand{\labelenumi}{\emph{(\roman{enumi})}}
\item \emph{Limit law.}
  At every continuity point of $\Psi_f^\Diag$, we have
\begin{equation}\label{LL002}
\Psi_f^\Diag(R)=\lim_{M\to\infty}\lambda\left\{X\in\T_\Theta^{\Diag} : M^{-n/2}|S_f(M,X)|>R\right\}.
\end{equation}

\item \emph{Tail asymptotic.}
We have
\begin{equation}\label{tailasymptotic002}
\Psi_f^\Diag(R)=C_{n,f}^\Diag R^{-4}(\log R)^{n-1}
+O_{n,f}\!\left(R^{-4}(\log R)^{n-2}\log \log R\right),
\qquad R\to\infty,
\end{equation}
where $C_{n,f}^{\Diag}$ is given in \eqref{eq:diagonal-tail-constant}, and in particular $C_{n,f}^{\Diag}$ is positive if and only if $P_\square f\neq0$.

\item \emph{Moments.}
If the density of $\lambda$ is in $L^\infty(\T_\Theta^\Diag)$, then for every $0<\beta<4$,
\begin{align}
  \lim_{M\to\infty}\int_{\T_\Theta^{\Diag}}\left(M^{-n/2}|S_f(M,X)|\right)^\beta\, \d\lambda(X)
  =\beta\int_0^\infty R^{\beta-1}\Psi_f^\Diag(R)\,\d R<\infty.
\label{eq:diagonal-moment-convergence-intro}
\end{align}
For $\beta=2$, the limit \eqref{eq:diagonal-moment-convergence-intro} evaluates to $2^n\|P_\square f\|_2^2$.
\end{enumerate}
\end{theorem}

Accounting for the different choices of normalisation, the tail asymptotic \eqref{tailasymptotic002} generalises the one obtained by Marklof in \cite{Marklof-spectral-form-factors} when $n = 2$.

\begin{remark}
 For the centred Gaussian $f(\vecx)=\e^{-\pi\|\vecx\|^2}$, the constants $C_{n,f}$ and $C_{n,f}^\Diag$ become
\begin{equation}\label{eq:gaussian-diagonal-constant-intro}
C_{n,f}= \frac{2}{(n+1)(2^n + 1)}\frac{\prod_{j=2}^n \xi(j)}{\prod_{j=1}^n\xi(2j)}, \qquad C_{n,f}^\Diag=\frac{2^{3n-2}}{\pi^n(n-1)!},
\end{equation}
where $\xi(s) := \pi^{-s/2}\Gamma(s/2) \zeta(s)$ denotes the completed Riemann zeta function.
\end{remark}

\begin{remark} 
The limits in \eqref{LL001} and \eqref{LL002} extend to bounded, rapidly decreasing, Riemann-integrable cut-offs; see Section \ref{sec:second}.  This includes
characteristic functions of bounded sets whose boundaries have Lebesgue
measure zero. We have not proved that the tail asymptotic
\eqref{tailasymptotic001} extends to more singular cut-off functions.
\end{remark}

Finally, fix $X\in\Sym_n(\R)$ and sample the one-parameter family $tX$.

\begin{theorem}\label{thm:main3} 
Let $n\geq2$, $f\in\Si(\R^n)$, and $X\in\Sym_n(\R)$ with
$\det X\neq0$.  There is a right-continuous non-increasing function
$\Psi_{f,X}:\R\to[0,1]$ with $\Psi_{f,X}(R)=1$ for $R<0$, such that, for every absolutely continuous
Borel probability measure $\lambda$ on $\R$,
\begin{equation}\label{LL003}
\Psi_{f,X}(R)=\lim_{M\to\infty}\lambda\left\{
t\in\R:M^{-n/2}|S_f(M,tX)|>R\right\}
\end{equation}
at every continuity point of $\Psi_{f,X}$.  Furthermore:
\begin{enumerate}
\renewcommand{\labelenumi}{\em (\roman{enumi})}
\item
If $\dim_{\Q}\operatorname{span}_{\Q}\{x_{ij}:1\leq i\leq j\leq n\}
>\frac{n(n-1)}2+1$, then
$\Psi_{f,X}=\Psi_f$. 
\item If $X=\diag(x_1,\ldots,x_n)$ with $x_i/x_j\notin\Q$ for every $i\neq j$, then $\Psi_{f,X}=\Psi_f^\Diag$.
\end{enumerate}
Consequently, the tails in the two cases are given by
\eqref{tailasymptotic001} and \eqref{tailasymptotic002}, respectively.
\end{theorem}

\begin{remark}
The hypothesis on $X$ in (i) is satisfied, for example, if its coefficients $x_{ij}$, where $i \leq j$, are linearly independent over $\Q$. In this case $\dim_{\Q}\operatorname{span}_{\Q}\{x_{ij}:1\leq i\leq j\leq n\} =\frac{n(n+1)}{2}$.
\end{remark}

\begin{remark}
Extending Theorem \ref{thm:main3} to the setting when $f$ is a sharp cut-off requires uniform second-moment estimates related to quantitative versions of the Oppenheim conjecture; see, for example, \cite{Eskin-Margulis-Mozes-II}. We do not cover this case here.
\end{remark}

\begin{remark}
When $X$ is positive definite and $f$ is a smooth approximation to an indicator, the normalised sum  $M^{-n}|S_f(M,tX)|^2$ is the corresponding normalised spectral form factor for the flat torus determined by $X$; cf.~\cite{Marklof-spectral-form-factors} for background.  Thus Theorem
\ref{thm:main3} gives a limit law for this smoothed form factor.
\end{remark}

The paper is organised as follows.  Section \ref{section:prelim} recalls the Shale--Weil representation, toward the construction of the theta function $\Theta_f$, and shows that its modulus is naturally defined on a homogeneous space of finite volume. The theta function is our main tool in relating the normalised theta sum to homogeneous dynamics.  Section \ref{section:closed-horospheres}  uses the equidistribution of horospheres to prove the limit laws in \eqref{LL001} and \eqref{LL002}. Using upper bounds for  the theta function we are then able to extend the limit laws to convergence of moments. Section \ref{section:one-parameter} extends the limit theorems of the previous sections to the setting of one-parameter families of quadratic forms.  Section \ref{sec:second} extends the full and diagonal convergence statements to Riemann-integrable cut-offs. Finally,  Sections \ref{section:Rank-n-asymptotics} and  \ref{section:diagonal-asymptotics} prove the tail asymptotics  \eqref{tailasymptotic001} and \eqref{tailasymptotic002}, respectively.

\subsection*{Acknowledgements}
We thank Nimish Shah for helpful discussions. This research was supported by EPSRC grants EP/S024948/1 and EP/W007010/1 and through access to OpenAI models provided by the ChatGPT for Academic Researchers program.

\subsection*{Data access statement}
No new data were generated or analysed during this study.

\subsection*{Declaration of AI use}

The authors produced the manuscript draft comprising the main results and proofs for this paper.  OpenAI's Codex was used to assist with mathematical exposition, drafting and clarification of proof arguments, and LaTeX editing. In particular, this allowed us to extend the calculation of the constant in the leading-order asymptotics \eqref{eq:generic-tail-constant} from $n=3$ to higher dimensions, and to relax the hypotheses on $X$ in Theorem \ref{cor:application-of-shah} (i). The authors take full responsibility for the work presented in this manuscript.

\section{Automorphic theta functions}\label{section:prelim}

We begin by recalling the background needed for the construction of the theta function $\Theta_f$ associated to a given Schwartz function $f\in\Si(\R^n)$. This is our main tool in the study of the value distribution of the original theta sum $S_f(M,X)$. For a detailed account  we refer the reader to  the classic monographs \cite{Lion-Vergne,Folland-Harmonic-Analysis}.

Throughout, we use row vectors and normalisation consistent with
\cite{Marklof-Welsh-higher-rank-theta-sums,Marklof-Welsh-higher-rank-theta-sumsII}.

We define the Heisenberg group $\Heisd$ as the set $\R^n \times \R^n \times \R$ with group multiplication given by 
\begin{equation}
    (\vecp, t)(\vecp', t') := (\vecp + \vecp', t + t' + \tha\vecp J \trans \vecp'), \qquad J := \matr{0}{-I}{I}{0}.\label{def:matrix-J}
\end{equation}
The symplectic group is defined by
\begin{equation}
    \spdr := \lcur g \in \mathrm{GL}(2n,\R) : g J \trans g = J \rcur.\label{def:spdr}
\end{equation}
For $g=\matr{A}{B}{C}{D}\in\spdr$, set $Z_0=\i I$.
We realise the \emph{metaplectic group} as pairs
\begin{equation}\label{def:metaplectic-group}
\mpdr=\left\{(g,q):g\in\spdr,\ q\in\mathrm U(1),\quad
q^2=\frac{\det(CZ_0+D)}{|\det(CZ_0+D)|}\right\}.
\end{equation}
The projection is $\mathrm{pr}(g,q)=g$; thus the two lifts of $g$
are $(g,q)$ and $(g,-q)$.

To define the multiplication, recall the Siegel upper half-space
\begin{equation}\label{def:Siegel-uhp}
\uhp_n:=\{X+\i Y:X,Y\in\Sym_n(\R),\ Y\text{ positive definite}\}
\end{equation}
and the action of $\spdr$ of $\uhp_n$ defined by
\begin{equation}\label{def:siegel-mobius}
\matr{A}{B}{C}{D}\cdot Z=(AZ+B)(CZ+D)^{-1}.
\end{equation}
For $Z\in\uhp_n$, the matrix $CZ+D$ can be shown to be invertible, and therefore this action is well defined. Let $j(g,Z):=\det(CZ+D)$. Since $\uhp_n$ is simply connected, each pair $(g,q)$ determines a unique holomorphic square root $\eta_{g,q}$ with
\begin{equation}\label{eq:metaplectic-square-root}
\eta_{g,q}(Z)^2=j(g,Z),\qquad
\eta_{g,q}(Z_0)=q\,|j(g,Z_0)|^{1/2}.
\end{equation}
The group law for $\mpdr$ in this parametrisation is then
\begin{equation}\label{eq:metaplectic-multiplication}
(g_1,q_1)(g_2,q_2)
=\left(g_1g_2,\;
q_2\frac{\eta_{g_1,q_1}(g_2\cdot Z_0)}
{|\eta_{g_1,q_1}(g_2\cdot Z_0)|}\right).
\end{equation}
Indeed, the corresponding square root is
$\eta_{g_1,q_1}(g_2\cdot Z)\eta_{g_2,q_2}(Z)$, by the identity
$j(g_1g_2,Z)=j(g_1,g_2\cdot Z)j(g_2,Z)$.
The identity is $(I,1)$, and the projection has central kernel
$\{(I,1),(I,-1)\}$.

We define a left action of $\spdr$ on $\Heisd$ by
\begin{equation}
    g \cdot (\vecp, t) := (\vecp g\inv, t),\label{def:spdr-action-on-H}
\end{equation}
where $g \in \spdr$, and $(\vecp, t) \in \Heisd$.
This action preserves $J$, and so is an action by automorphisms.
 The Jacobi group is then defined as the semidirect product
\begin{align}
     \Jac = \Heisd \rtimes \spdr,\label{def:jacobi}
\end{align}
with multiplication given by 
\begin{align}\label{def:jacobi-multiplication}
    (h, g)(h', g') = (h (g \cdot h'), g g'),
\end{align}
where $g, g' \in \spdr$, $h, h' \in \Heisd$, and $g \cdot h'$ is given by \eqref{def:spdr-action-on-H}.
The \emph{metaplectic Jacobi group} is $\Heisd\rtimes\mpdr$,
where $(g,q)$ acts on $\Heisd$ through $g$ and
\eqref{def:spdr-action-on-H}.

Let $((\vecx, \vecy), t) \in \Heisd$. The Schr\"odinger representation $W$ of $\Heisd$ on $\LtRd$ by unitary operators is defined by
\begin{equation}
    W(((\vecx, \vecy), t)) f(\vecw) = e\!\lp -t + \tha \vecx \trans\vecy + \vecw\trans\vecy\rp f(\vecw + \vecx).\label{def:schrod-repn}
\end{equation}
For $g\in\spdr$, let $W^g(h):=W(g^{-1}\cdot h)$.  By the Stone--von Neumann theorem there is a unitary operator $R(g)$ such that
\begin{equation}
W^g(h)=R(g)^{-1}W(h)R(g),
\qquad
W(g\cdot h)=R(g)W(h)R(g)^{-1}.
\label{def:proj-schrod-weil}
\end{equation}
Consequently, there is a unitary cocycle $\rho : \spdr \times \spdr \to S^1$ such that
\begin{equation}
    R(g g') = \rho (g, g')R(g) R(g').\label{eq:cocycle}
\end{equation}

The projective \emph{Schr\"odinger--Shale--Weil representation} of the Jacobi group $\Jac$ is defined by $R(h,g):=W(h)R(g)$. Indeed,
\begin{equation}\label{eq:jacobi-projective-law}
\begin{aligned}
R((h,g)(h',g')) & =W(h)W(g\cdot h')R(g g') \\
& = \rho(g,g') W(h) R(g) W(h') R(g') \\
& = \rho(g,g') R(h,g) R(h',g').
\end{aligned}
\end{equation}

For $f \in \Si(\R^n)$ we formally define the theta function  $\Theta_f : \Jac \to \C$ by the infinite sum
\begin{align}
    \Theta_f (h, g) = \sum_{\vecm \in \Z^n} [R(h,g)f](\vecm)\label{def:theta-fn}.
\end{align}
We see in  Proposition~\ref{prop:theta-well-defined} that for $f \in \Si(\R^n)$ this series converges pointwise.

Let $Y$ be a positive definite matrix, and let $U$ be an upper triangular, unipotent matrix such that
\begin{equation}
    Y = U V \trans U \label{def:Y-to-V},
\end{equation}
where $\boldsymbol v=(v_1,\ldots,v_n)\in\R_{>0}^n$ and
$V=\diag(\boldsymbol v)$.  Let $\boldsymbol v^{1/2}=(v_1^{1/2},\ldots,v_n^{1/2})$ and
\mbox{$V^{1/2}:=\diag(\boldsymbol v^{1/2})$}.  We define
\begin{equation}\label{def:Y-half}
    Y^{1/2} := U V^{1/2},
\end{equation}
so that $Y^{1/2} \trans Y^{1/2} = Y$, and define
$Y^{-1/2}:=(Y^{1/2})^{-1}$.
Let $K := \spdr \cap \mathrm{O}(2n)$. We have that $K \simeq \Un$ via 
 \begin{align}
    Q \mapsto k(Q) := \matr{\Re (Q)}{-\Im (Q)}{\Im (Q)}{\Re (Q)}. \label{eq:K-to-Un}
 \end{align}
The Iwasawa decomposition of $\spdr$ with respect to $K$ is given by
\begin{align}
    g = \matr{I}{X}{0}{I} \matr{Y^{\tha}}{0}{0}{\trans Y^{-\tha}}k(Q)\label{def:Iwasawa-decomp}
\end{align}
where $X,Y\in\Sym_n(\R)$, $Y$ is positive definite, $Q \in \Un$, and $k(Q)$ as in \eqref{eq:K-to-Un}. It will sometimes be convenient to use the notation
\begin{equation}
  \label{eq:n(X)-a(Y)}
 n(X) := \matr{I}{X}{0}{I}, \qquad a(Y) := \matr{Y^{\tha}}{0}{0}{\trans Y^{-\tha}}.
\end{equation}

For a parabolic element $p = \matr{A}{B}{0}{\trans A\inv} \in \spdr$, we have by \cite[Proposition 2.1]{Marklof-Welsh-higher-rank-theta-sums} that
\begin{align}\label{SLaction}
   R \lp \matr{A}{B}{0}{\trans A\inv}\rp f (\vecw) = |\det A|^{\tha} e\!\lp \tha \vecw A \trans B \trans \vecw \rp f(\vecw A).
 \end{align}

 In this normalisation we have that if $p=\smatr{A}{B}{0}{\trans A\inv}$, then $\rho(p,g)=\rho(g,p)=1$, and so the cocycle is trivial when either argument is parabolic. See \cite[Proposition 2.1]{Marklof-Welsh-higher-rank-theta-sums}.

Using this, we may write the theta function in Iwasawa coordinates as
\begin{multline}
    \Theta_f(((\vecx, \vecy), t), \tg) = (\det Y)^{1/4}e\!\lp -t + \tha \vecx \trans\vecy\rp \sum_{\vecm \in \Z^n} f_{Q}((\vecm + \vecx)Y^{1/2})\\ \times e\!\lp\tha (\vecm + \vecx)X \trans(\vecm + \vecx) + \vecm \trans \vecy\rp \label{def:theta-fn-iwasawa}
  \end{multline}
where $f_Q := R(k(Q))f$. Note that the modulus $|\Theta_f|$ and more generally the product $\Theta_{f_1}\overline{\Theta_{f_2}}$ are independent of $\rho$. Equation~\eqref{def:theta-fn-iwasawa} yields the connection to the Siegel theta sum, which we summarise in the following proposition.

\begin{proposition}\label{prop:relation-to-theta-sums}
    Let $M>0$, $X\in\Sym_n(\R)$, and let $f:\R^n\to\C$ satisfy $|f(\vecw)|\leq C (1 + \|\vecw\|)^{-(n+1)}$ for some $C<\infty$. Then
        \begin{align}\label{eq:theta-sum-theta-relation}
        M^{-n/2}S_f(M,X) = \Theta_f(\vecnull,g_{X, M})
    \end{align}
    with
    \begin{equation}\label{eq:theta-sum-Iwasawa-point}
       g_{X, M}=\matr{I}{X}{0}{I} \matr{M\inv I}{0}{0}{M I}.
    \end{equation}
  \end{proposition}

For $g=n(X)a(Y)k(Q)$ as in \eqref{def:Iwasawa-decomp}, its lower
blocks satisfy
\begin{equation}
CZ_0+D=\trans Y^{-1/2}Q,\qquad
\frac{j(g,Z_0)}{|j(g,Z_0)|}=\det Q.
\end{equation}
Thus the pair $(g,q)\in\mpdr$ has Iwasawa coordinates
$(X+\i Y,Q,q)$ with $q^2=\det Q$.
In particular, the inverse image $\KMp=\mathrm{pr}^{-1}(K)$ is
identified, via $(k(Q),q)\leftrightarrow(Q,q)$, with
\begin{equation}\label{eq:metaplectic-compact}
\begin{gathered}
\KMp\simeq\{(Q,q):Q\in\Un,\ q\in\mathrm U(1),\ q^2=\det Q\},\\
(Q,q)(Q',q')=(QQ',qq').
\end{gathered}
\end{equation}
Here the product follows from \eqref{eq:metaplectic-multiplication} since every element of $K$ fixes $Z_0$. For parabolic elements $p = \smatr{A}{B}{0}{\trans A\inv}$ with $\det A > 0$ we always take the lift $(p,1)$, which is continuous in $p$.

We distinguish the projective operators $R(g)$ from the genuine metaplectic operators $\RMp(g,q)$. One explicit phase convention
is obtained from the Gaussian $h_0(\vecw)=\e^{-\pi\|\vecw\|^2}$. Since $R(k(Q))h_0=\alpha(Q)h_0$ with $|\alpha(Q)|=1$, we may take
\begin{equation}\label{eq:metaplectic-phase}
\begin{gathered}
\alpha(Q)=(R(k(Q))h_0)(\bn),\qquad
\chi(g,q)=q^{-1}\alpha(Q)^{-1},\\
\RMp(g,q)=\chi(g,q)R(g).
\end{gathered}
\end{equation}
With this choice $\RMp(k(Q),q)h_0 = q\inv h_0$ for all $(Q, q) \in \KMp$.
This fixes the usual metaplectic representation, and defines
\begin{equation}\label{eq:lifted-theta}
\ThetaMp_f(h,(g,q))
:=\sum_{\vecm\in\Z^n}[W(h)\RMp(g,q)f](\vecm)
=\chi(g,q)\Theta_f(h,g).
\end{equation}
We point out that as the phase is independent of $f$, it does not appear in the formulae for $|\Theta_f|$ and $\Theta_{f_1}\overline{\Theta_{f_2}}$. The explicit phase is given by \eqref{eq:diagonal-metaplectic-phase}.

\begin{proposition}\label{prop:theta-well-defined}
    Let $f \in \Si (\R^n)$. Then the series \eqref{def:theta-fn} is absolutely convergent. Furthermore, $\ThetaMp_f$ is smooth on $\Heisd\rtimes\mpdr$. The modulus $|\Theta_f|$ and, for $f_1,f_2\in\Si(\R^n)$, the mixed products $\Theta_{f_1}\overline{\Theta_{f_2}}$ are continuous on $\Jac$.
\end{proposition}    
    
\begin{proof}   
    By \eqref{def:theta-fn-iwasawa}, we have that
    \begin{align}
        |\Theta_f(h, \tg)| &\leq (\det Y)^{1/4} \sum_{\vecm \in \Z^n} |f_Q ((\vecm + \vecx)Y^{1/2})| \notag\\
        &\ll_{f,A} (\det Y)^{1/4} \sum_{\vecm \in \Z^n} \left(1+\|(\vecm +\vecx)Y^{1/2}\|\right)^{-A},\label{eq:theta-abs-conv}
    \end{align}
where the final estimate follows from the fact that if $f \in \Si (\R^n)$, then so is $f_Q$, and in particular for any $A > 0$ and any multi-index $\alpha$ we have \cite[Lemma 4.2]{Marklof-Welsh-higher-rank-theta-sums}
    \begin{align}\label{eq:uniform-schwartz-seminorm}
      \labs \lp\frac{\partial}{\partial \vecw}\rp^\alpha f_Q (\vecw) \rabs \leq c_f(\alpha, A) (1 + \|\vecw\|)^{-A},
    \end{align}
    uniformly in $Q$. Taking $A>n$ shows that the series defining $\ThetaMp_f$ is absolutely convergent. Since $\Mp (n, \R)$ acts smoothly on $\Si (\R^n)$ via $\RMp$, we may differentiate term by term, and the differentiated series converges locally uniformly. It follows that $\ThetaMp_f$ is smooth. Since $|\Theta_f|$ and $\Theta_{f_1}\overline{\Theta_{f_2}}$ are independent of the phase in \eqref{eq:lifted-theta}, they are well-defined continuous functions on $\Jac$.
\end{proof}

It is clear from \eqref{def:theta-fn-iwasawa} that the cocycle $\rho$ and the Heisenberg variable $t$ play no role in the value distribution of $|\Theta_f|$. Therefore, we can view $|\Theta_f|$ as a function on
\begin{equation}
 \AffSym:=\R^{2n}\rtimes\spdr 
\end{equation}
    Here $\spdr$ acts on $\R^{2n}$ from the left by 
    \begin{equation}
        g \cdot \vecp := \vecp g\inv, \qquad g \in \spdr, \quad \vecp \in \R^{2n},\label{def:spdr-r2n-action}
    \end{equation}
     which leads to the multiplication
    \begin{align}
        (\vecp, g) (\vecp' , g') = (\vecp + g\cdot \vecp' , g g').\label{def:affsym-mult} 
    \end{align}
      We identify $\spdr$ with
$\uhp_n\times\Un$ via the Iwasawa decomposition, where $\uhp_n$ and its $\spdr$-action are given
in \eqref{def:Siegel-uhp} and \eqref{def:siegel-mobius}.

Let $\gamma=\smatr{A}{B}{C}{D}\in\spdz$, and write $h_\gamma=(\vecp_\gamma,0)$, where $\vecp_\gamma=(\vecr,\vecs)\in\{0,\tha\}^{2n}$ is determined by
\begin{equation}\label{eq:theta-characteristic}
  2r_i\equiv(C\trans D)_{ii} \pmod2, \qquad 2s_i \equiv(A\trans B)_{ii}\pmod2.
\end{equation}
For $u=((\vecm,\vecn),t)$, $\vecm,\vecn\in\Z^n$, the exact transformation law shown in \cite[Theorem 4.1]{Marklof-Welsh-higher-rank-theta-sums} is
\begin{equation}\label{eq:exact-theta-automorphy}
\Theta_f((u h_\gamma,\gamma)(h,g))
=\varepsilon(\gamma)\rho(\gamma,g)
e\!\lp-t+\tha\vecm\trans\vecn\rp\Theta_f(h,g),
\qquad |\varepsilon(\gamma)|=1.
\end{equation}
Motivated by this, we define the following  subgroup of $\AffSym$,
\begin{align}
\Gamma_\Theta = \{ ((\vecm, \vecn) + \vecp_{\gamma}, \gamma) : \gamma \in \spdz, \ \vecm, \vecn \in \Z^{n}\}.\label{eq:Gamma-explicit-definition}
\end{align}

In view of \eqref{eq:exact-theta-automorphy}, we have the following proposition.

\begin{prop}\label{prop:modulus-of-theta}
Let $f_1,f_2\in\Si(\R^n)$. Then for every $\gamma \in \Gamma_\Theta$,
\begin{equation}\label{eq:theta-product-invariance}
\Theta_{f_1}(\gamma\cdot(\vecp,g))
\overline{\Theta_{f_2}(\gamma\cdot(\vecp,g))}
=\Theta_{f_1}(\vecp,g)\overline{\Theta_{f_2}(\vecp,g)}
\end{equation}
\end{prop}

We embed $\spdr$ in $\AffSym$ via the map
\begin{equation}\label{eq:symplectic-affine-embedding}
\iota: \spdr\hookrightarrow \AffSym, \qquad g \mapsto (\vecnull, g) .
\end{equation}

The following lemma shows that we can also view the modulus of $\Theta_f(g):=\Theta_f(\vecnull, g)$ as a function on $\spdr$ which is invariant under a lattice $\overline\Gamma_\Theta<\spdr$.

\begin{lemma}\label{Lem1}
The group 
\begin{align}
    \overline\Gamma_\Theta = \{ \gamma \in \spdz : \vecp_{\gamma}=\vecnull\}
        \label{eq:Gamma'-explicit-definition}
\end{align}
has finite index in $\spdz$.
\end{lemma}

\begin{proof}
We have the explicit representation
\begin{align}
    \overline\Gamma_\Theta 
    =\left\{ \matr{A}{B}{C}{D} \in \spdz : 
    (A\trans B)_{ii} \equiv (C\trans D)_{ii} \equiv 0 \bmod 2,  1\leq i\leq n \right\}.
    \label{eq:Gamma'-explicit-definition2}
\end{align}

From this it is clear that $\oGam$ has finite index in $\spdz$, since it contains the principal congruence subgroup of elements $\gamma\equiv I \pmod 2$.
\end{proof}

Observe that  the group
\begin{equation}\label{eq:theta-lattice-intersection}
\Gamma_\Theta\cap \iota\spdr
\end{equation}
is a lattice in $\iota\spdr$ and we have the identifications
\begin{equation}\label{eq:theta-orbit-quotient}
\Gamma_\Theta\bs \Gamma_\Theta\, \iota\spdr 
\simeq (\Gamma_\Theta\cap \iota\spdr)\bs \iota\spdr
\simeq  \overline\Gamma_\Theta\bs \spdr.
\end{equation}

Thus, we can represent the space $\Gamma_\Theta\bs \Gamma_\Theta\, \iota\spdr$ as a fibre bundle over the quotient $\spdz\bs\spdr$, where each fibre is given by the finite set $\overline\Gamma_\Theta\bs \spdz$. A fundamental domain for $\Gamma_\Theta\bs \Gamma_\Theta\, \iota\spdr$ (i.e., a fundamental domain for the $\Gamma_\Theta$ action on the set $\Gamma_\Theta\, \iota\spdr$) is given by
\begin{equation}\label{FDTheta}
\widehat\calG_\Theta = \{ (\vecp_{\gamma^{-1}},g) : \gamma\in  \Lambda_\Theta, \; g\in\calG_{\spdz} \} \subset \Gamma_\Theta\, \iota\spdr,
\end{equation}
where $\calG_{\spdz}$ is a fundamental domain for $\spdz$ in $\spdr$, and $\Lambda_\Theta$ is a set of representatives of the finite set $\overline\Gamma_\Theta\bs \spdz$.

\section{Limit theorems and convergence of moments}
\label{section:closed-horospheres}
\label{section:limit-theorems}
\label{sec:moment-convergence}

The main limit theorem for theta sums $S_f(M, X)$ when $f \in \Si (\R^n)$ essentially follows from the equidistribution of closed horospheres. However, as the theta function is unbounded, we cannot prove convergence of moments directly in this way.

We use the height function on $\spdz \bs \spdr$, pulled back to the quotient $\oGam \bs \spdr$, defined by
\begin{equation}\label{eq:modular-height}
H(\oGam g):=\max_{\gamma\in\spdz}
  (\det Y(\gamma g))^{1/4}.
\end{equation}
By \cite[Corollary 4.5]{Marklof-Welsh-higher-rank-theta-sums},
\begin{equation}\label{eq:theta-estimate}
|\Theta_f(\oGam g)|\ll_{n,f}H(\oGam g).
\end{equation}

Throughout we let $\mu_\Theta$ denote the $\spdr$-invariant probability measure on $\oGam \bs \spdr$. We have the following estimate.
\begin{lemma}[\cite{Marklof-Welsh-higher-rank-theta-sums}, Lemma 5.1]
  For $0 < \beta < 2(n + 1)$,
\begin{equation}
  \label{eq:tail-decay-D}
  \mu_\Theta( \oGam g \in \oGam \bs \spdr : H^\beta(\oGam g) > R) \ll R^{-\tfrac{2 (n + 1)}{\beta}}.
\end{equation}
In particular, $H^\beta \in L^1 (\oGam \bs \spdr)$ for $\beta$ in the given range.
\end{lemma}

Recall $n(X)$ from \eqref{eq:n(X)-a(Y)}, and define
\begin{equation}
  \label{eq:diagonal-flow}
  \Phi_t := \matr{\e^{t/2}I}{0}{0}{\e^{-t/2}I}.
\end{equation}
It is classical that the flow defined by $\Phi_t$ is mixing on $\oGam \bs \spdr$, with respect to $\mu_\Theta$.

\begin{theorem}\label{thm:hgt-fn-integral} For $0 < \beta < 2(n+1)$, and $\lambda \ll m_\Theta$ with a density $\rho \in L^\infty(\T_\Theta)$, we have 
  \begin{equation}
    \label{eq:hgt-fn-estimate}
   \sup_{M \geq 1} \int_{\T_\Theta} H\lp\oGam n(X) \Phi_{-2\log M}\rp^\beta \: \de \lambda (X) < \infty.
  \end{equation}
\end{theorem}

\begin{proof}
Since $\d\lambda\leq\|\rho\|_\infty\d m_\Theta$, it suffices to prove the theorem for $m_\Theta$. Let
\begin{equation}
P^-:=\left\{\matr{A}{0}{0}{\trans A\inv}
\matr{I}{0}{T}{I}:A\in\GL(n,\R),\ T\in\Sym_n(\R)\right\},
\end{equation}
and let $B_\delta\subset P^-$ be a sufficiently small neighbourhood
defined by $\|A-I\|<\delta$ and $\|T\|<\delta$.
By \cite[Lemma 5.3]{Marklof-Welsh-higher-rank-theta-sums},
$H(\oGam gh)\asymp_n H(\oGam g)$ uniformly for $h\in B_\delta$.
Moreover,
\begin{equation}
\Phi_{2\log M}B_\delta\Phi_{-2\log M}\subseteq B_\delta,
\qquad M\geq1,
\end{equation}
because conjugation replaces $T$ by $M^{-2}T$ and fixes $A$.
Consequently,
\begin{align}
&\int_{\T_\Theta}H(\oGam n(X)\Phi_{-2\log M})^\beta\,\d m_\Theta(X) \notag \\
&\quad\ll_{n,\beta}\int_{\T_\Theta}\int_{B_\delta}
H(\oGam n(X)h\Phi_{-2\log M})^\beta\,\d h\,\d m_\Theta(X).
\end{align}
The map $(X,h)\mapsto\oGam n(X)h$ is locally a diffeomorphism.
Since the closed horosphere and $\overline{B_\delta}$ are compact,
the push-forward of the measure in the last integral has a bounded
density with respect to $\mu_\Theta$. This also follows from the Haar
coordinate formula in \cite[Section 2.1]{Marklof-Welsh-higher-rank-theta-sums}. Therefore
\begin{equation}
  \int_{\T_\Theta}\int_{B_\delta} H(\oGam n(X)h\Phi_{-2\log M})^\beta\,\d h\,\d m_\Theta(X)
  \ll  \int_{\oGam \bs \spdr} H(\oGam g \Phi_{-2\log M})^\beta\de \mu_\Theta(\oGam g). 
\end{equation}
Right invariance of Haar measure therefore gives
\begin{equation}
  \int_{\oGam \bs \spdr}   H(\oGam g \Phi_{-2\log M})^\beta\de \mu_\Theta(\oGam g) = \int_{\oGam\bs\spdr}H^\beta\,\d\mu_\Theta<\infty,
\end{equation}
uniformly for $M\geq1$, as required.
\end{proof}

The following theorem is the equidistribution of horospheres stated in our precise case. It follows from mixing and Margulis' thickening argument, using the same transverse neighbourhood as in the preceding proof. For a more sophisticated approach, including rates of convergence, see  \cite{Edwards-effective-equidist} and references therein.

\begin{theorem}\label{thm:equidist-of-horospheres}
Let $n\geq 1$ and let $\lambda \ll m_\Theta$ be a Borel probability measure on $\T_\Theta$.  For every $F\in C_b(\overline\Gamma_\Theta\bs\spdr)$,
\begin{equation}\label{eq:full-horosphere-equidistribution}
\lim_{M\to \infty}\int_{\T_\Theta}F(\overline\Gamma_\Theta n(X)\Phi_{-2 \log M})\,\d\lambda(X) =\int_{\overline\Gamma_\Theta\bs\spdr}F\,\d\mu_\Theta.
\end{equation}
\end{theorem}

Theorem \ref{thm:hgt-fn-integral} and Theorem \ref{thm:equidist-of-horospheres} along with a classical argument from probability theory give the main theorem of this section.

\begin{theorem}\label{thm:full-moment-convergence}
Let $n\geq1$, $f\in\Si(\R^n)$, and let $\lambda\ll  m_\Theta$ be a
Borel probability measure on $\T_\Theta$. For any $F\in C_b(\R_{\geq0})$,
\begin{equation}\label{eq:full-limit-law}
\lim_{M\to\infty}\int_{\T_\Theta}F\!\left(M^{-n/2}|S_f(M,X)|\right)\,\d\lambda(X)
=\int_{\overline\Gamma_\Theta\bs\spdr}
F\circ|\Theta_f|\,\d\mu_\Theta.
\end{equation}
If, in addition, $\lambda$ has a density $\rho \in L^\infty(\T_\Theta)$, then, for every
$0<\beta<2(n+1)$,
\begin{equation}\label{eq:full-moment-convergence}
\lim_{M\to\infty}\int_{\T_\Theta}\left(M^{-n/2}|S_f(M,X)|\right)^\beta\,\d\lambda(X)=\int_{\overline\Gamma_\Theta\bs\spdr}|\Theta_f|^\beta\,\d\mu_\Theta.
\end{equation}
\end{theorem}

\begin{proof}
  The limit given in equation~\eqref{eq:full-limit-law} follows immediately from Theorem~\ref{thm:equidist-of-horospheres} and Proposition~\ref{prop:relation-to-theta-sums}.

    For the second claim, we first define the set
  \begin{equation}
    \label{eq:tail-set}
    \mathcal{T}(R) := \{\oGam g \in \oGam \bs \spdr : H^\beta(\oGam g) > R\},
  \end{equation}
  and let $\chi_R$ be a smoothed indicator of $\mathcal T(R + 1)$, satisfying
  \begin{align}
    \label{eq:smoothe-tail-cut-off}
    &\chi_R(\oGam g) =
    \begin{dcases*}
      1 & if $\oGam g \in \mathcal{T}(R+1)$\\ 0 & if $\oGam g\not\in \mathcal{T}(R)$,
    \end{dcases*}\\
     &\text{and } 0 \leq \chi_R(\oGam g) \leq 1 \text{ if } \oGam g \in \mathcal{T}(R) \bs \mathcal{T}(R + 1)
  \end{align}
  Then,
  \begin{multline}
    \int_{\T_\Theta} \big|\Theta_f (\oGam n(X) \Phi_{-2\log M})\big|^\beta \de \lambda (X) =
    \int_{\T_\Theta} [|\Theta_f|^\beta\cdot \chi_R] (\oGam n(X) \Phi_{-2\log M}) \de \lambda (X)\\
    + \int_{\T_\Theta} [|\Theta_f|^\beta\cdot(1 - \chi_R)] (\oGam n(X) \Phi_{-2\log M}) \de \lambda (X).
  \end{multline}
  As $|\Theta_f|^\beta (1 - \chi_R)$ is bounded, we apply Theorem~\ref{thm:equidist-of-horospheres} to see that
  \begin{equation}
    \lim_{ M \to \infty}\int_{\T_\Theta} [|\Theta_f|^\beta\cdot(1 - \chi_R)] (\oGam n(X) \Phi_{-2\log M}) \de \lambda (X) = \int_{\oGam \bs \spdr} |\Theta_f|^\beta\cdot(1 - \chi_R) \de \mu_\Theta.
  \end{equation}
    Also, it is clear from \eqref{eq:tail-decay-D} that
  \begin{equation}
    \lim_{R \to \infty}\int_{\oGam \bs \spdr} |\Theta_f|^\beta \chi_R \: \de \mu_\Theta =  0.
  \end{equation}
  Writing $\beta = 2(n + 1) - s$, we have by \eqref{eq:theta-estimate} that 
  \begin{equation}
    |\Theta_f|^\beta \chi_R \ll H^{2(n + 1) - s}\chi_R = H^{2(n + 1) - \tfrac{s}{2}}\chi_R H^{-\tfrac{s}{2}} \leq H^{2(n + 1) - \tfrac{s}{2}}R^{-\tfrac{s}{2\beta}}.
  \end{equation}
  Now we apply Theorem \ref{thm:hgt-fn-integral} to see that
  \begin{equation}
    \sup_{M \geq 1} \,\int_{\T_\Theta} [|\Theta_f|^{\beta} \cdot \chi_R](\oGam n(X) \Phi_{-2 \log M}) \: \de \lambda (X) \ll R^{-\frac{s}{2\beta}}.
  \end{equation}
  Therefore, given $\veps > 0$, we may take $R$ and then $M$ sufficiently large so that
  \begin{equation}
     \bigg|\:\int_{\T_\Theta} |\Theta_f(\oGam n(X) \Phi_{-2\log M})|^\beta \de \lambda (X) - \int_{\oGam \bs \spdr} |\Theta_f|^\beta \de \mu_\Theta \:\bigg| < \veps.
  \end{equation}
\end{proof}

When $\beta = 2$ we have more precise information.

\begin{prop}\label{prop:covariance-limit}
  If $f_1,f_2:\R^n\to\C$ are bounded, rapidly decreasing, and Riemann-integrable on every compact rectangle, then
  \begin{equation}
    \label{eq:theta-sum-variance}
    \lim_{M \to \infty}\frac{1}{M^n} \int_{\T_\Theta}S_{f_1}(M,X) \overline{S_{f_2}(M,X)} \, \de m_\Theta (X) = \int_{\R^n} f_1(\vecx) \overline{f_2 (\vecx)} \, \de \vecx +  \int_{\R^n} f_1(\vecx) \overline{f_2 (-\vecx)} \, \de \vecx.
  \end{equation}
\end{prop}

\begin{proof}
 By our choice of cut-off functions $f_1, f_2$, the sums $S_{f_1}(M, X)$ and $S_{f_2}(M,X)$ are defined by absolutely convergent series. Therefore we have that
\begin{multline}\label{eq:sum-of-integrals}
    \int_{\T_\Theta} S_{f_1}(M,X) \overline{S_{f_2}(M,X)} \de m_\Theta(X) 
    \\=\sum_{\vecm_1, \vecm_2 \in \Z^n} f_1 (\tfrac1M \vecm_1 )\overline{f_2 (\tfrac1M\vecm_2)} \int_{\T_\Theta} e(\tha(\vecm_1 X \trans \vecm_1 - \vecm_2 X \trans \vecm_2)) \de m_\Theta(X).
  \end{multline}
We note that
\begin{multline}
  \int_{\T_\Theta} e\big(\tha(\vecm_1 X \trans \vecm_1 - \vecm_2 X \trans \vecm_2)\big) \de m_\Theta(X) \\
    = \frac{1}{2^n}\lp\prod_{i = 1}^n \int_0^2 e ( \tha x_{ii} (m_{1,i}^2 - m_{2,i}^2)) \de x_{ii}\rp \lp\prod_{i< j} \int_0^1 e ( x_{ij} (m_{1,i}m_{1,j} - m_{2,i} m_{2,j})) \de x_{ij}\rp.
\end{multline}
Therefore
\begin{equation}
\int_{\T_\Theta}e\!\left(\tfrac12(\vecm_1 X\trans\vecm_1 - \vecm_2 X\trans\vecm_2)\right)\,\d m_\Theta(X)
=\ind_{\{\trans\vecm_1\vecm_1=\trans\vecm_2\vecm_2\}}.
\end{equation}
Note that we have $\trans\vecm_1\vecm_1=\trans\vecm_2\vecm_2$ if and only if $\vecm_1=\vecm_2$ or $\vecm_1=-\vecm_2$. It follows that
\begin{align}\label{eq:sum-of-integrals2}
&M^{-n}\int_{\T_\Theta}S_{f_1}(M,X)\overline{S_{f_2}(M,X)}\,\d m_\Theta(X)
\notag\\
&=M^{-n}\lp\sum_{\vecm\in\Z^n}f_1(M^{-1}\vecm)  \left(\overline{f_2(M^{-1}\vecm)}+\overline{f_2(-M^{-1}\vecm)}\right)\rp
  - M^{-n}f_1(\vecnull)\overline{f_2(\vecnull)}.
\end{align}
The last term tends to zero as $M \to \infty$, and the remaining expression is a Riemann sum. Taking $M \to \infty$  gives \eqref{eq:theta-sum-variance}.

\end{proof}

Recall that when defining  $\Theta_f$ as a function on a quotient, we must fix a metaplectic lift, and use the same choice of lift for every $f$. That is, for $\mathscr X_\Theta=\oGam\bs\spdr$, we define $s_\Theta^{\Mp}:\mathscr X_\Theta\to\mpdr$ as a Borel measurable map that selects a representative of each coset together with one of its two possible lifts, and define
\begin{equation}\label{eq:scalar-theta-section}
\Theta_f(x):=\ThetaMp_f(\vecnull,s_\Theta^{\Mp}(x)),
\qquad x\in\mathscr X_\Theta.
\end{equation}
Changing the lift multiplies all these functions by the same phase factor, as discussed in \eqref{eq:exact-theta-automorphy}.
Therefore the modulus, mixed products and $L^2$ identities are independent of the choices.

\begin{prop}\label{L2}
With the convention \eqref{eq:scalar-theta-section}, the map $f\mapsto\Theta_f$, initially defined for $f\in\Si(\R^n)$, extends uniquely to a bounded linear map $L^2(\R^n)\to L^2(\mathscr X_\Theta,\mu_\Theta)$.  Furthermore for $f_1,f_2\in L^2(\R^n)$, we have that 
  \begin{equation}
    \label{eq:Theta-L2-identity}
    \int_{\oGam \bs \spdr} \Theta_{f_1}\overline{\Theta_{f_2}} \, \de \mu_\Theta = \int_{\R^n}f_1(\vecw) \overline{f_2 (\vecw)} \, \de \vecw +  \int_{\R^n}f_1(\vecw) \overline{f_2 (-\vecw)} \, \de \vecw.
  \end{equation}
\end{prop}

\begin{proof}
  First assume $f_1, f_2 \in \Si (\R^n)$. Apply Theorem~\ref{thm:full-moment-convergence} with $\beta=2$ to $f=f_1+\i^k f_2$, $k=0,1,2,3$, and use the complex polarisation identity. We obtain
  \begin{equation}
    \int_{\oGam \bs \spdr} \Theta_{f_1}\overline{\Theta_{f_2}} \, \de \mu_\Theta = \lim_{M \to \infty}\frac{1}{M^n} \int_{\T_\Theta}S_{f_1}(M,X) \overline{S_{f_2}(M,X)} \, \de m_\Theta (X).  
   \end{equation}
   Proposition \ref{prop:covariance-limit} shows  \eqref{eq:Theta-L2-identity} in this setting. We then extend by density to complete the proof.
 \end{proof}

 \begin{theorem}\label{thm:theta-sum-variance} Let $n \geq 1$ and let $\lambda \ll m_\Theta$ be a Borel probability measure with density $\rho \in L^\infty (\T_\Theta)$. Suppose $f_1, f_2 \in \mathcal{S}(\R^n)$. Then 
   \begin{equation}
     \lim_{M \to \infty}\frac{1}{M^n} \int_{\T_\Theta}S_{f_1}(M,X) \overline{S_{f_2}(M,X)} \, \de \lambda (X) = \int_{\R^n}f_1(\vecw) \overline{f_2 (\vecw)} \, \de \vecw +  \int_{\R^n}f_1(\vecw) \overline{f_2 (-\vecw)} \, \de \vecw.
   \end{equation}
 \end{theorem}

 \begin{proof}
   By complex polarisation of Theorem~\ref{thm:full-moment-convergence} at $\beta=2$, applied to $f=f_1+\i^k f_2$ for $k=0,1,2,3$, we have
   \begin{equation}
     \lim_{M \to \infty}\frac{1}{M^n} \int_{\T_\Theta}S_{f_1}(M,X) \overline{S_{f_2}(M,X)} \, \de \lambda (X) = \int_{\oGam \bs \spdr} \Theta_{f_1}\overline{\Theta_{f_2}} \, \de \mu_\Theta.
   \end{equation}
   Therefore the previous proposition gives the result.
 \end{proof}

 Now suppose
\begin{equation}
  X := \diag(\vecx) =
  \begin{pmatrix}
    x_1 & & \\ & \ddots & \\ & & x_n
  \end{pmatrix}, 
\end{equation}
where $\vecx := ( x_1, \dots, x_n) \in \R^n$.

Recall that in \eqref{def:theta-fn} we define the theta function as
\begin{align}
    \Theta_f (h,g) = \sum_{\vecm \in \Z^n} [W(h)R(\tg)f](\vecm),
\end{align}
where $W$ is the Schr\"odinger representation of the Heisenberg group $\Heisd$, and $R$ is the projective Shale--Weil representation of the symplectic group $\spdr$. We may study $\Theta_f$ restricted to a subgroup of $\Jac$ by simply restricting the Schr\"odinger--Shale--Weil representation to the given subgroup.

We embed $(\sltr)^n$ into $\spdr$ via the map
\begin{align}\label{eq:sltr-to-spdr}
   \psi: \left(\matr{a_1}{b_1}{c_{1}}{d_{1}}, \dots, \matr{a_{n}}{b_{n}}{c_{n}}{d_{n}} \right) \longmapsto \matr{A}{B}{C}{D}
\end{align}
where $\smatr{a_i}{b_i}{c_i}{d_i} \in \sltr$, and 
\begin{align}
    &A := \diag (a_1, \dots, a_n), ~~
    B := \diag (b_1, \dots, b_n),\notag \\
    &C := \diag (c_1, \dots, c_n), ~~
    D:= \diag (d_1, \dots, d_n).
\end{align}
We then define the projective Shale--Weil representation on $\prodsltr$ as the projective Shale--Weil representation $R$ on $\spdr$ restricted to $\psi(\prodsltr)$.

Let $g \in \sltr$. Then there exist $u \in \R$, $v > 0$ and $\phi \in [0,2\pi)$ such that
\begin{align}
    g = n_u a_v k_{\phi},
 \end{align}
where
 \begin{align}
    n_u = \matr{1}{u}{0}{1}, a_v = \matr{v^{1/2}}{0}{0}{v^{-1/2}}, k_\phi := \matr{\cos \phi}{ -\sin \phi}{\sin \phi }{\cos \phi}.
 \end{align}
Let $g_1, \dots, g_n \in \sltr$, and let
\begin{align}
    g_i = n_{u_i} a_{v_i}k_{\phi_i}
\end{align}
be the Iwasawa decomposition of $g_i$, $1 \leq i \leq n$. Let 
\begin{align}
    \vecu := (u_1, \dots, u_n),\qquad \vecv:= (v_1,\dots, v_n), \qquad \vecphi := (\phi_1, \dots, \phi_n),
\end{align}
and 
\begin{align}
    n(\vecu) := \psi(n_{u_1}, \dots, n_{u_n}), a(\vecv):= \psi(a_{v_1},\dots, a_{v_n}), k(\vecphi) := \psi(k_{\phi_1}, \dots, k_{\phi_n}). \label{def:iwasawa-vectors}
\end{align}
Using that $R$ is defined uniquely up to a phase, the restriction of the projective representation to $\prodsltr$ can then be described explicitly through the following identities from the projective Shale--Weil representation on $\sltr$:
\begin{align}
    [R(n(\vecu))f](\vecw) &:= e(\tha \vecw \:\diag(\vecu) \trans\vecw) f(\vecw),\label{eq:schrod-weil-diag-hor}\\
    [R(a(\vecv))f](\vecw) &:= |\det \diag (\vecv)|^{1/4}f(\vecw\diag(\vecv)^{1/2}),\label{eq:schrod-weil-diag-geo}\\
    [R(k(\vecphi))f](\vecw) &:= \int_{\R^n} G_{\phi_1}(w_1, w'_1) \cdots G_{\phi_n}(w_n, w'_n) f(\vecw') \de \vecw'. \label{eq:schrod-weil-diag-rot}
\end{align}
The integral formula is valid when every $\sin \phi_i\neq0$. If $\phi_i=0\pmod{2\pi}$, the operator is the identity in that coordinate;
if $\phi_i=\pi\pmod{2\pi}$, it is reflection in that coordinate.
For $\sin \phi\neq0$, the kernel is
\begin{equation}
    G_{\phi}(w,w'):= \frac{1}{\sqrt{|\sin \phi|}}e\!\lp\frac{\tha(w^2 + w'^2)\cos \phi - ww'}{\sin \phi}\rp. \label{eq:Greens-fn1}
 \end{equation}
For real lifted angles $\phi$, define
\begin{equation}\label{eq:diagonal-metaplectic-phase}
\sigma(\phi)=
\begin{cases}
2\nu,&\phi=\nu\pi,\\
2\nu+1,&\nu\pi<\phi<(\nu+1)\pi,
\end{cases}
\quad \nu\in\Z,
\end{equation}
and set $q(\vecphi):=\e^{\frac{\i}{2}\sum_{j=1}^n\phi_j}$.
The pair $(k(\vecphi),q(\vecphi))\in\mpdr$ then satisfies
\begin{equation}
\RMp(k(\vecphi),q(\vecphi))
=e\!\left(-\frac18\sum_{j=1}^n\sigma(\phi_j)\right)R(k(\vecphi)).
\end{equation}
In one coordinate this is the phase $e(-\sigma(\phi)/8)$ that appears in
\cite[Section 2.1]{Cellarosi-Marklof}.  At $\phi=\nu\pi$, the lifted operator is
$e(-\nu/4)f((-1)^\nu w)$. Thus the lift is smooth across multiples
of $\pi$, and changes sign under $\phi\mapsto\phi+2\pi$.
In formulae involving $|\Theta_f|$, and $\Theta_{f_1}\overline{\Theta_{f_2}}$ we write $f_{\vecphi}:=R(k(\vecphi))f$, as this phase plays no role. Equivalently, let
$Q(\vecphi)=\diag(\e^{\i\phi_1},\ldots,\e^{\i\phi_n})$.
Then $k(\vecphi)=k(Q(\vecphi))$ and $f_{\vecphi}=f_{Q(\vecphi)}$.

Let $f \in \Si (\R^n)$, and suppose $(g_1, \dots, g_n) \in \prodsltr$, with Iwasawa decompositions $g_i = n_{u_i}a_{v_i}k_{\phi_i}$.  Set $z_i = u_i + \i v_i$, and
\begin{align}\label{eq:diag-coordinates}
    \vecz = (z_1, \dots, z_n), \: \vecu = (u_1, \dots, u_n), \: \vecv = (v_1, \dots, v_n), \: \vecphi = ( \phi_1, \dots, \phi_n).
\end{align}
Using \eqref{eq:schrod-weil-diag-hor}--\eqref{eq:schrod-weil-diag-rot}, we may write the restriction of $|\Theta_f|$ to $\AffSL$ in `Iwasawa coordinates' as a non-negative, real-valued function on $\R^n \times \R^n \times (\uhp_1\times [0,2\pi))^n$, given by
\begin{multline}
    |\Theta_f ((\vecx, \vecy); \vecu + \i\vecv, \vecphi)| :=\\ (\det \diag (\vecv))^{\frac{1}{4}} \labs\sum_{\vecm \in \Z^n} \!\!f_{\vecphi}((\vecm + \vecx) (\diag (\vecv))^{\ha}) e\!\lp\frac{1}{2}(\vecm + \vecx)\diag(\vecu)\trans(\vecm + \vecx) + \vecm\trans\vecy\rp\rabs,\label{def:theta-diag-iwasawa}
\end{multline}
where $f_{\vecphi}:=R(k(\vecphi))f$ and $\diag(\vecv)^{1/2}:=\diag(v_1^{1/2},\ldots,v_n^{1/2})$.

By the general invariance properties of  $|\Theta_f|$ as a function on $\spdr$, we have that $|\Theta_f|$ restricted to $\AffSL$ is invariant under the action of
\begin{equation}
    \Delta_\Theta := \Gamma_\Theta \cap (\AffSL)
\end{equation}
where $\Gamma_\Theta$ is as defined in \eqref{eq:Gamma-explicit-definition}. We note that if $\pi : \AffSL \to \sltrn$ denotes the standard projection, then $\pi(\Delta_\Theta) = \sltzn$. We define
\begin{equation}
  \label{eq:oDel}
  \oDel := \oGam \cap (\sltr)^n,
\end{equation}
and set $\nu_\Theta$ to be the $\sltrn$-invariant probability measure on the quotient $\oDel \bs \sltrn$.

 Note $\{\Phi_t : t \in \R\} < \sltrn < \spdr$, and further, the right action of $\Phi_t$ is mixing on $\oDel \bs \sltrn$ with respect to $\nu_\Theta$. We now state the analogue of Theorem \ref{thm:equidist-of-horospheres} in this setting.

\begin{theorem}\label{thm:equidist-of-horospheres-diag}
Let $n\geq1$ and let $\lambda\ll m_\Theta^{\Diag}$ be a Borel probability
measure on $\T_\Theta^{\Diag}$.  For every $F\in C_b(\oDel \bs (\sltr)^n)$,
\begin{equation}
  \lim_{t\to -\infty}\int_{\T_\Theta^{\Diag}}F\!\left(\oDel n(\vecx) \Phi_{t}\right)\,\d\lambda(\vecx) =\int_{\oDel \bs (\sltr)^n} F\,\d\dlimmeas.
\label{eq:diagonal-horosphere-equidistribution}
\end{equation}
\end{theorem}

This gives the limit theorem for diagonal $X$. For moments, let
\begin{equation}
H_1(g):=\max_{\gamma\in\sltz}\Im(\gamma g\cdot\i)^{1/4},
\qquad D(g_1,\ldots,g_n):=\prod_{j=1}^n H_1(g_j).
\end{equation}
Note that $H_1$ is $\sltz$-invariant in each factor, and so $D$ is defined on the quotient $\oDel\bs\sltrn$. By \eqref{def:theta-diag-iwasawa} we have
\begin{equation}
  |\Theta_f ((\vecx, \vecy); \vecu + \i\vecv, \vecphi)| \leq
  (\det \diag (\vecv))^{\frac{1}{4}} \sum_{\vecm \in \Z^n} \!\!\labs f_{\vecphi}((\vecm + \vecx) (\diag (\vecv))^{\ha})\rabs,
\end{equation}
and so if  $f \in \Si (\R^n)$, then it is not hard to see that
\begin{equation}
  \label{eq:4}
  |\Theta_f|\ll_{n,f}D.
\end{equation}
Clearly  $\oDel$ is the product of $n$ copies of the theta subgroup
\begin{equation}
  \Big\{g \in \sltz : g \equiv \matr{1}{0}{0}{1} \text{ or } \matr{0}{1}{1}{0} \pmod 2\Big\} < \sltz,
\end{equation}
so $\nu_\Theta$ is a product measure. By Lemma~\ref{lem:simplex} we have that
\begin{equation}\label{eq:diag-height-estimate}
\nu_\Theta\{D>R\}\ll_n R^{-4}(\log R)^{n-1},\qquad R\geq2.
\end{equation}
In particular, $D^\beta\in L^1(\nu_\Theta)$ for $0<\beta<4$. The estimate \eqref{eq:diag-height-estimate} is sharp. Indeed, for $n\geq2$ choose a Schwartz function with $P_\square f\neq0$. The inclusion $\{|\Theta_f|>B_fR\}\subseteq\{D>R\}$, with $|\Theta_f|\leq B_fD$, and Theorem~\ref{thm:diagonal-tail} give the matching lower bound $\nu_\Theta\{D>R\}\gg_{n,f}R^{-4}(\log R)^{n-1}$. For $n=1$ this follows directly from integrating $\de u \de v / v^2$ over the cusp $v > R^4$.

Repeating the proof of Theorem~\ref{thm:hgt-fn-integral} in this setting gives the following.

\begin{prop}\label{thm:hgt-fn-integral-diag} For $0 < \beta < 4$, and $\lambda \ll m_\Theta^{\Diag}$ with  density $\rho \in L^\infty(\T_\Theta^{\Diag})$, we have 
  \begin{equation}
    \label{eq:hgt-fn-estimate-diag}
   \sup_{M \geq 1} \int_{\T_\Theta^{\Diag}} D\lp\oDel n(X) \Phi_{-2\log M}\rp^\beta \: \de \lambda (X) < \infty.
  \end{equation}
\end{prop}

Replacing Theorem~\ref{thm:equidist-of-horospheres} with Theorem~\ref{thm:equidist-of-horospheres-diag}, and Theorem~\ref{thm:hgt-fn-integral} with Proposition~\ref{thm:hgt-fn-integral-diag} leads to the convergence of moments in the context of diagonal forms.

\begin{theorem}\label{thm:diagonal-moment-convergence}
Let $n\geq1$, $f\in\Si(\R^n)$, and let $\lambda\ll m_\Theta^{\Diag}$ be a Borel probability measure on $\T_\Theta^{\Diag}$.  For every $F\in C_b(\R_{\geq0})$,
\begin{equation}
  \lim_{M\to\infty}\int_{\T_\Theta^{\Diag}} F\!\left(M^{-n/2}|S_f(M,X)|\right)\,\d\lambda(X)
  =\int_{\oDel \bs(\sltr)^n} F(|\Theta_f|)\, \d\dlimmeas.
\label{eq:diagonal-limit-law}
\end{equation}
If $\lambda$ has density $\rho\in L^\infty(\T_\Theta^{\Diag})$, then, for every $0<\beta<4$,
\begin{equation}\label{eq:diagonal-moment-convergence}
  \lim_{M\to\infty}\int_{\T_\Theta^{\Diag}} \left(M^{-n/2}|S_f(M,X)|\right)^\beta\, \d\lambda(X)
  =\int_{\oDel \bs (\sltr)^n}|\Theta_f|^\beta\, \d\dlimmeas.
\end{equation}
\end{theorem}

\begin{prop}\label{cor:L2-product}
  The map $f\mapsto\Theta_f$ extends to a bounded linear map $L^2(\R^n)\to L^2(\oDel \bs \sltrn,\nu_\Theta)$. For $f_1,f_2\in L^2(\R^n)$, we have
\begin{equation}\label{eq:L2-product}
\begin{aligned}
\int_{\oDel \bs (\sltr)^n} \Theta_{f_1}\overline{\Theta_{f_2}}\,\d\dlimmeas =
\sum_{\boldsymbol\varepsilon\in\{\pm1\}^n}\int_{\R^n}f_1(\vecw)\overline{f_2(\boldsymbol\varepsilon\vecw)}\,\d\vecw.
\end{aligned}
\end{equation}
\end{prop}

\begin{proof}
  First let  $f_1,f_2\in\Si(\R^n)$, and consider the integral
  \begin{multline}
    M^{-n}\int_{\T_\Theta^\Diag}  S_{f_1}(M,X)\overline{S_{f_2}(M,X)}\,\d m_\Theta^\Diag(X) \\= \frac{1}{M^{n}}\sum_{\vecm_1, \vecm_2 \in \Z^n} f_1 (\tfrac1M \vecm_1 )\overline{f_2 (\tfrac1M\vecm_2)} \int_{\T^\Diag_\Theta}e(\tha(\vecm_1 X \trans \vecm_1 - \vecm_2 X \trans \vecm_2)) \de m_\Theta^{\Diag}(X).
  \end{multline}
We note 
\begin{equation}
  \int_{\T^\Diag_\Theta} e\big(\tha(\vecm_1 X \trans \vecm_1 - \vecm_2 X \trans \vecm_2)\big) \de m_\Theta^\Diag(X)
  = \frac{1}{2^n}\prod_{i = 1}^n \int_0^2 e\big( \tha x_{ii} (m_{1,i}^2 - m_{2,i}^2)\big) \de x_{ii}. 
\end{equation}
This product is non-zero only when  $m_{1, i} = \pm m_{2,i}$ for each $i$, and therefore, in a similar way to \eqref{eq:sum-of-integrals2} we have 
\begin{multline}
    M^{-n}\int_{\T_\Theta^\Diag}  S_{f_1}(M,X)\overline{S_{f_2}(M,X)}\,\d m_\Theta^\Diag(X) \\= \frac{1}{M^{n}}\sum_{\boldsymbol\varepsilon \in \{\pm 1\}^n}\sum_{\vecm \in \Z^n} f_1 (\tfrac1M \vecm )\overline{f_2 (\tfrac1M\boldsymbol\varepsilon\vecm)} + O(M^{-1}). 
  \end{multline}
Therefore,
  \begin{equation}\label{eq:diag-covariance-limit}
\lim_{M \to \infty}M^{-n}\int_{\T_\Theta^\Diag}
 S_{f_1}(M,X)\overline{S_{f_2}(M,X)}\,\d m_\Theta^\Diag(X) =\sum_{\boldsymbol\varepsilon\in\{\pm1\}^n}
 \int_{\R^n}f_1(\vecw)\overline{f_2(\boldsymbol\varepsilon\vecw)}\,\d\vecw.
\end{equation}
The same torus calculation is valid for bounded, rapidly decreasing,
Riemann-integrable cut-offs. For Schwartz functions, complex polarisation of Theorem~\ref{thm:diagonal-moment-convergence} at $\beta=2$ gives \eqref{eq:L2-product}. As  $\Si (\R^n)$ is dense in $L^2(\R^n)$, and  $\|\Theta_f\|_2^2=2^n\|P_\square f\|_2^2$ for $f \in \Si(\R^n)$, we may extend to all of $L^2(\R^n)$.
\end{proof}

This leads to the analogue of Theorem \ref{thm:theta-sum-variance}.

\begin{theorem}\label{thm:diagonal-covariance} Let $n \geq 1$ and let $\lambda \ll m^\Diag_\Theta$ be a Borel probability measure with density $\rho \in L^\infty (\T^\Diag_\Theta)$. Suppose $f_1, f_2 \in \mathcal{S}(\R^n)$. Then 
   \begin{equation}
     \lim_{M \to \infty}\frac{1}{M^n} \int_{\T^\Diag_\Theta}S_{f_1}(M,X) \overline{S_{f_2}(M,X)} \, \de \lambda (X) =  \sum_{\boldsymbol{\veps}\in\{\pm 1\}^n} \int_{\R^n} f_1(\vecw) \overline{f_2(\boldsymbol{\veps}\vecw)} \, \de \vecw.
   \end{equation}
 \end{theorem}

\begin{proof}
Apply complex polarisation to Theorem~\ref{thm:diagonal-moment-convergence}
at $\beta=2$, and then use \eqref{eq:L2-product}.
\end{proof}

\section{Limit theorems for one-parameter families}
\label{section:one-parameter}
To study the distribution of $\frac{1}{M^{n/2}}|S_f(M, tX)|$, for fixed $X \in \Sym_n(\R)$, with $\det X\neq 0$ and random $t$, we replace the equidistribution of closed horospheres with the equidistribution of expanding translates of unipotent orbit segments. This rests on Ratner's orbit closure theorems \cite[Theorems~A and~B]{Rat91} and Shah’s equidistribution theorem \cite[Theorem 1.4]{Shah-expanding-translates}. We state these results in the context of $\spdr$.

\begin{theorem}[Ratner, Shah]\label{thm:application-of-shah}
Let $\Gamma$ be a lattice in $\spdr$, and let $\varphi:\sltr\hookrightarrow\spdr$ be a continuous embedding.
There is a closed connected subgroup $H<\spdr$ containing $\varphi(\sltr)$ such that
\begin{enumerate}
\renewcommand{\labelenumi}{\em (\roman{enumi})}
\item $\Gamma_H:=\Gamma \cap H$ is a lattice in $H$;
\item $\Gamma_H\bs H \simeq \Gamma \bs \Gamma H$ carries the
$H$-invariant probability measure $\mu_H$;
\item $\overline{\Gamma \bs \Gamma \varphi(\sltr)} = \Gamma \bs \Gamma H$;
\item 
for $F \in C_b (\Gamma\bs \spdr)$ and $\lambda$ an absolutely continuous Borel probability measure  on $\R$, we have
    \begin{align}\label{eq:shah-one-parameter}
        \lim_{y \to 0}\int_{\R} F\lp \Gamma\varphi\lp\matr{1}{t}{0}{1}\matr{y^{1/2}}{0}{0}{y^{-1/2}}\rp\rp \de \lambda(t)
        = \int_{\Gamma_H\bs H} F(\Gamma h)\,\d\mu_H(h).
    \end{align}
\end{enumerate}
\end{theorem}

\begin{proof}
Ratner's orbit-closure theorem gives (i)--(iii), since $\varphi(\sltr)$ is generated by one-parameter unipotent subgroups. Applying \cite[Theorem 1.4]{Shah-expanding-translates} on $\Gamma_H\bs H$, where the image of $\varphi(\sltr)$ is dense, gives (iv). \end{proof}

 Let
    \begin{equation}\label{eq:sl2-triple}
        \frakn_+(B) = \matr{0}{B}{0}{0}, ~~~ \frakn_{-}(C) = \matr{0}{0}{C}{0}, ~~~ \fraka (T):= \matr{T}{0}{0}{-\trans T}.
    \end{equation}
   The matrices $\frakn_+(X), \frakn_-(X\inv)$, and $\fraka (I)$ form an $\mathfrak{sl}_2$-triple, i.e., satisfy the commutation relations
   \begin{equation}\label{eq:sl2-commutators}
[\fraka(I),\frakn_+(X)] = 2 \frakn_+(X), \qquad 
[\fraka (I),\frakn_-(X\inv)] = -2 \frakn_-(X\inv), \qquad 
[\frakn_+(X),\frakn_-(X\inv)] = \fraka(I) .
\end{equation}
Exponentiating these yields an embedding $\varphi_X : \sltr \hookrightarrow \spdr$ given by 
\begin{equation}
\varphi_X\!\left(\matr{a}{b}{c}{d}\right)
=\matr{aI}{bX}{cX^{-1}}{dI}.
\end{equation}
In particular, we have
    \begin{equation}\label{eq:sl2-embedded-flow}
        \varphi_X \lp \matr{1}{t}{0}{1}\matr{y^{1/2} }{0}{0}{y^{-1/2}}\rp = \matr{I}{tX}{0}{I}\matr{y^{1/2} I}{0}{0}{y^{-1/2}I} .
    \end{equation}
Let
\begin{equation}\label{eq:one-parameter-orbit-closure}
\mathcal O_X:= \overline{\Gamma\bs\Gamma\varphi_X(\sltr)}.
\end{equation}
By Theorem~\ref{thm:application-of-shah} we have that $\mathcal O_X=\Gamma\bs\Gamma H$, with $m_X$ the $H$-invariant probability measure on $\mathcal{O}_X$. Theorem \ref{cor:application-of-shah} below gives necessary and sufficient conditions to determine when $H = \spdr$ or $\sltrn$. This reduces to determining the Lie algebra $\mathfrak{h}$ of $H$. As $\varphi_X(\sltr) \leq  H$, we have that $\mathfrak{a}(I)\in\mathfrak{h}$. Therefore $\mathfrak h$ is invariant under $\mathrm{ad}\: \mathfrak a (I)$. Furthermore
\begin{equation}
[\fraka(I),\frakn_+(B)] = 2 \frakn_+(B), \qquad 
[\fraka (I),\frakn_-(C)] = -2 \frakn_-(C), \qquad 
[\fraka (I),\fraka(T)] = 0.
\end{equation}
i.e. the eigenspaces of  $\mathrm{ad}\: \mathfrak a (I)$ in $\mathfrak{sp}(n,\R)$ are $ \frakn_+(\Sym_n(\R))$, $ \frakn_-(\Sym_n(\R))$ and $ \fraka (\mathrm{M}_n (\R))$, and so
\begin{equation}
 \mathfrak h =  (\frakn_+(\Sym_n(\R)) \cap \mathfrak h) \oplus (\frakn_-(\Sym_n(\R)) \cap \mathfrak h) \oplus (\fraka(\mathrm{M}_n(\R)) \cap \mathfrak h).
\end{equation}
We also observe that
\begin{equation}\label{eq:doublebracket}
[[\frakn_+(B),\frakn_-(X^{-1})],\frakn_+(C)] =\frakn_+(BX^{-1}C+CX^{-1}B),
\end{equation}
 and so  $\{B:\frakn_+(B)\in\mathfrak{h}\}$ is closed under the product 
\begin{equation}
  \label{eq:X-product}
  B \circ_X C := \frac12(BX^{-1}C+CX^{-1}B),
\end{equation}
and therefore is a Jordan algebra with unit $X$.

Conversely, the Kantor-Koecher-Tits construction turns a Jordan algebra $\mathcal{J}$ with unit $X$ for this product into the graded Lie subalgebra
\begin{equation}
  \label{eq:graded-X-lie-alg}
  \frakn_+(\mathcal{J})\oplus
  [\frakn_+(\mathcal{J}),\frakn_-(X^{-1}\mathcal{J}X^{-1})]\oplus
  \frakn_-(X^{-1}\mathcal{J}X^{-1})
\end{equation}
of $\mathfrak{sp}(n,\mathbb{R})$, which contains the triple \eqref{eq:sl2-commutators}. Given $X$, the relevant Jordan algebra $\mathcal{J} = \mathcal{J}_\Q(X)$ is defined as follows. Let $r(X):=\dim_{\Q}\operatorname{span}_{\Q}\{x_{ij}:i\leq j\}$ and write 
\begin{equation}\label{eq:rational-coefficient-space}
X=\sum_{\nu=1}^{r(X)}\xi_\nu A_\nu,
\end{equation}
where the $\xi_\nu$ form a $\Q$-basis of the span of the entries of $X$
and $A_\nu\in\Sym_n(\Q)$. Thus  
\begin{equation}
\mathcal V_{\Q}(X):=\operatorname{span}_{\Q}\{A_1,\ldots,A_{r(X)}\}
\end{equation}
is the smallest rational linear subspace whose real extension contains $X$. 

The polynomial $\det(\sum_\nu t_\nu A_\nu)$ has rational coefficients and is nonzero at $(\xi_1,\ldots,\xi_{r(X)})$. Hence it is nonzero at some rational tuple, so $\mathcal V_{\Q}(X)$ contains an invertible matrix. Choose any such $U$ and define $\mathcal J_{\Q}(X)$ as the Jordan algebra generated by $\mathcal V_{\Q}(X)$, with the Jordan product $\circ_U$ defined by
\begin{equation}\label{eq:JQ-operation}
B\circ_U C:=\frac12(BU^{-1}C+CU^{-1}B).
\end{equation}
Note that $U$ is the unit of $\mathcal J_{\Q}(X)$. 

  The Jordan algebra $\mathcal{J}_{\Q}(X)$ does not depend on the choice of $U$. Indeed, let $U'\in\mathcal{V}_{\mathbb{Q}}(X)$ be another invertible matrix. The space $\mathcal{J}_{\mathbb{Q}}(X)U^{-1}$ is closed under the product $\frac12(ab+ba)$. As it contains $I$ and all powers of $U'U^{-1}$, it follows by the Cayley--Hamilton theorem that $(U'U^{-1})^{-1} \in \mathcal{J}_\Q(X)U\inv$, and so $UU'^{-1}U\in\mathcal{J}_{\mathbb{Q}}(X)$. The algebra $\mathcal{J}_{\mathbb{Q}}(X)$ is then closed under $\circ_{U'}$, and so contains the Jordan algebra generated by $\mathcal{V}_{\mathbb{Q}}(X)$ under $\circ_{U'}$, and by symmetry the two coincide.

\begin{theorem}\label{cor:application-of-shah}
Let $\Gamma$ be commensurable with $\Sp(n,\Z)$, and let $X\in\Sym_n(\R)$ have $\det X\neq0$.  Then, for any $F\in C_b(\Gamma\bs\spdr)$ and any absolutely continuous Borel
probability measure $\lambda$ on $\R$,
    \begin{align}\label{eq:one-parameter-horosphere}
        \lim_{y \to 0}\int_{\R} F\lp \Gamma\matr{I}{tX}{0}{I}\matr{y^{1/2} I}{0}{0}{y^{-1/2}I}\rp \de \lambda(t) = \int_{\Gamma\bs \spdr} F \de m_X .
    \end{align}
In particular:
\begin{enumerate}
\renewcommand{\labelenumi}{\em (\roman{enumi})}
\item 
The measure $m_X$ is the $\spdr$-invariant probability measure on
$\Gamma\bs\spdr$ if and only if
\begin{equation}\label{eq:jordan-density-criterion}
\mathcal J_{\Q}(X)=\Sym_n(\Q).
\end{equation}
In particular, \eqref{eq:jordan-density-criterion} holds if
\begin{equation}\label{eq:rational-rank-condition}
r(X)>\frac{n(n-1)}2+1.
\end{equation}

\item
For $X=\diag(\boldsymbol x)$, the measure $m_X$ is the $\psi(\sltr^n)$-invariant probability measure on $\Gamma\bs\Gamma\psi(\sltr^n)$ if and only if 
\begin{equation} \label{eq:diagonal-ratio-condition}
\frac{x_i}{x_j}\notin\Q \quad \text{for every $i\neq j$,} 
\end{equation}
where $\psi$ is as defined in \eqref{eq:sltr-to-spdr}.
\end{enumerate}
\end{theorem}

Case (ii) is also discussed in \cite[Lemma 3.2]{Kim-Marklof-Welsh-BT}.

\begin{proof}

\textbf{Case (i).} Let $\mathbf H$ be the smallest algebraic
$\Q$-subgroup of $\mathrm{Sp}_n$ whose real points contain $\varphi_X(\sltr)$. In view of \cite[Proposition 3.2 and Remark 3.7]{Shah-uniform-orbits}, we have that $H=\mathbf H(\R)^\circ$, the connected component of $\mathbf H(\R)$ containing the identity. In particular, $\mathfrak h=\operatorname{Lie}(H)$ is defined over $\Q$.

Let
\begin{equation}
\mathcal L_{\Q}:=\{B\in\Sym_n(\Q):\frakn_+(B)\in\mathfrak h\} .
\end{equation}
Since $\frakn_+(X)\in\mathfrak h$, and $\mathfrak h$ is defined over $\Q$, we have that $\mathcal V_{\Q}(X)\subseteq\mathcal L_{\Q}$. As $\frakn_-(X\inv) \in \mathfrak{h}$, then \eqref{eq:doublebracket} shows that
\begin{equation}
  \mathcal L_{\R}:=\mathcal L_{\Q}\otimes_{\Q}\R  =\{B\in\mathrm{Sym}_n(\mathbb{R}):\mathfrak{n}_+(B)\in\mathfrak{h}\}
\end{equation}
is a Jordan algebra with respect to the product $\circ_X$, with unit $X$.
The Weyl element
\begin{equation}
w_X=\varphi_X\!\lp\matr{0}{1}{-1}{0}\rp
=\matr{0}{X}{-X^{-1}}{0}
\end{equation}
satisfies
\begin{equation}\label{eq:w-x-conjugation}
  w_X\,\mathfrak{n}_+(B)\,w_X^{-1}=-\mathfrak{n}_-(X^{-1}BX^{-1}).
\end{equation}

Choose an invertible $U \in\mathcal{V}_\Q (X) \subseteq \mathcal{L}_\R$. We have, after arguing as in the paragraph after \eqref{eq:JQ-operation}, that $XU^{-1}X\in\mathcal L_{\R}$. It then follows from \eqref{eq:w-x-conjugation} that $-\mathfrak{n}_-(U^{-1})\in\mathfrak{h}$. With $U^{-1}$ in place of $X^{-1}$ in \eqref{eq:doublebracket}, and using that $U$ is rational, we find that $\mathcal L_{\Q}$ is closed under $\circ_U$. Thus $\mathcal J_{\Q}(X)\subseteq\mathcal L_{\Q}$.

If  $\mathcal{J}_\Q (X) = \Sym_n(\Q)$ then $\mathfrak{n}_+(\mathrm{Sym}_n(\mathbb{R}))\subseteq\mathfrak{h}$. Conjugating by $w_X$ gives $\mathfrak{n}_-(\mathrm{Sym}_n(\mathbb{R}))\subseteq\mathfrak{h}$, and the brackets $[\mathfrak{n}_+(B),\mathfrak{n}_-(C)]=\mathfrak{a}(BC)$ span $\mathfrak{a}(M_n(\mathbb{R}))$. Hence $\mathfrak{h}=\mathfrak{sp}(n,\mathbb{R})$, so $H =\mathrm{Sp}_n(\R)$ and $m_X$ is the invariant probability measure on $\Gamma\backslash\mathrm{Sp}(n,\mathbb{R})$.

Now suppose $\mathcal{J}:=\mathcal J_{\Q}(X)\neq \Sym_n(\Q)$. Following the  Kantor-Koecher-Tits construction we define
\begin{equation}
\mathcal{K}=U^{-1}\mathcal{J}U^{-1},\qquad
\mathcal{D}=\operatorname{span}_{\Q}\{BC:B\in \mathcal{J},\ C\in \mathcal{K}\},
\end{equation}
and
\begin{equation}\label{eq:jordan-lie-obstruction}
\mathfrak l_{\mathcal{J}}=\frakn_+(\mathcal{J})\oplus \fraka(\mathcal{D})\oplus \frakn_-(\mathcal{K}).
\end{equation}
For $B,C, E\in\mathcal{J}$ we define
\begin{equation}\label{eq:triple}
  \{B, C, E\}_U := \tfrac12(BU^{-1}CU^{-1}E+EU^{-1}CU^{-1}B).
\end{equation}
A calculation gives that
\begin{equation}
  \{B, C, E\}_U = (B\circ_U C)\circ_U E+(E\circ_U C)\circ_U B -(B\circ_U E)\circ_U C\in \mathcal{J}.
\end{equation}
Let $C=U^{-1}C'U^{-1}$, $F=U^{-1}F'U^{-1}$ in $\mathcal{K}$. Then
\begin{align}
  &[\mathfrak{a}(BC),\mathfrak{n}_+(E)]=2\,\mathfrak{n}_+\bigl(\{B,C',E\}_U\bigr),\\
  &[\mathfrak{a}(BC),\mathfrak{n}_-(F)]
  =-2\,\mathfrak{n}_-\bigl(U^{-1}\{C',B,F'\}_UU^{-1}\bigr),
\end{align}
so the bracket with $\mathfrak{a}(\mathcal{D})$ preserves
$\mathfrak{n}_+(\mathcal{J})$ and $\mathfrak{n}_-(\mathcal{K})$. Moreover
$[\mathfrak{n}_+(\mathcal{J}),\mathfrak{n}_-(\mathcal{K})]=\mathfrak{a}(\mathcal{D})$
by definition, the spaces $\mathfrak{n}_\pm(\cdot)$ are abelian, and
$[\mathfrak{a}(\mathcal{D}),\mathfrak{a}(\mathcal{D})]\subseteq\mathfrak{a}(\mathcal{D})$, and therefore $\mathfrak{l}_{\mathcal{J}}$ is a Lie subalgebra of $\mathfrak{sp}(n,\Q)$.

We observe that $X\in\mathcal{J}\otimes\mathbb{R}$, $U\in\mathcal{J}$, and  $U^{-1}\in\mathcal{K}$, with $[\mathfrak{n}_+(U),\mathfrak{n}_-(U^{-1})]=\mathfrak{a}(I)$. Arguing as in the paragraph after \eqref{eq:JQ-operation} shows that $X^{-1}\in\mathcal{K}\otimes\mathbb{R}$, and so it follows that  $\mathfrak{l}_{\mathcal{J}}\otimes_{\mathbb{Q}}\mathbb{R}$ contains the triple \eqref{eq:sl2-commutators}.

As $\varphi_X(\mathrm{SL}(2,\mathbb{R}))$ is generated by $\exp(\mathbb{R}\,\mathfrak{n}_+(X))$ and $\exp(\mathbb{R}\,\mathfrak{n}_-(X^{-1}))$, the adjoint action of $\varphi_X(\mathrm{SL}(2,\mathbb{R}))$ preserves $\mathfrak{l}_{\mathcal{J}}\otimes\mathbb{R}$. That is, $\varphi_X(\mathrm{SL}(2,\mathbb{R}))$ lies in the real points of the normaliser $\mathbf{N}$ of $\mathfrak{l}_{\mathcal{J}}$ in $\mathrm{Sp}_n$, which is an algebraic $\mathbb{Q}$-subgroup. This normaliser is proper, hence $\mathbf H$ is a proper subgroup of $\mathrm{Sp}_n$, and $\mathcal{O}_X \neq \Gamma \bs \spdr$.

For \eqref{eq:rational-rank-condition}, we note that a proper Jordan subalgebra of complex symmetric $n\times n$ matrices has codimension at least $n-1$ by \cite[Theorem 4.1]{Bik-Eisenmann-Sturmfels}. After complexification and a congruence taking $U$ to the identity, this applies to $\mathcal{J}$ and gives $\dim_{\Q}\mathcal{J}\leq \frac{n(n+1)}{2}-n+1$. Since $\dim_{\Q}\mathcal{J}\geq r(X)$, the stated inequality forces $\mathcal{J}=\Sym_n(\Q)$.

\textbf{Case (ii).} 
Let
\begin{equation}
D = (d_1,\ldots,d_n), \qquad d_j = \begin{pmatrix} x_j & \\ & 1 \end{pmatrix} . 
\end{equation}
Recall that $\psi$ denotes the diagonal embedding of $\sltr$ in $\spdr$, given in \eqref{eq:sltr-to-spdr}. We have
\begin{equation}
\varphi_X(\sltr) = \psi(\{ (d_1 g d_1^{-1},\ldots,d_n g d_n^{-1}) : g\in \sltr \}) ,
\end{equation}
which is the $D$ conjugate of the diagonally embedded $\sltr$ in $\sltr^n$.

The connected subgroups of $\sltr^n$ containing the diagonally embedded $\sltr$ are of the form 
\begin{equation}
H_{\mathcal P}=\{ (g_1,\ldots,g_n) \in \sltr^n : \text{$g_i=g_j$ if $i,j\in B_r$ for some $r$}\}
\end{equation}
where ${\mathcal P}=\{ B_1,\ldots, B_k\}$ is any partition of $\{1,\ldots, n\}$. Therefore, the relevant orbit is
\begin{equation}
 \Gamma\bs\Gamma \psi( \tilde H_{\mathcal P}), \qquad \tilde H_{\mathcal P}:=D H_{\mathcal P} D^{-1}, 
\end{equation}
for some partition $\mathcal P$. The orbit $ \Gamma\bs\Gamma \psi( \tilde H_{\mathcal P})$ is closed if and only if $\psi( \tilde H_{\mathcal P})\cap\Gamma$ is a lattice in $\psi( \tilde H_{\mathcal P})$. This is the same as saying (by the assumed commensurability) that $\tilde H_{\mathcal P} \cap\sltz^n$ is a lattice in $\tilde H_{\mathcal P}$, and hence that $H_{\mathcal P} \cap D^{-1} \sltz^n D$ is a lattice in $H_{\mathcal P}$. This holds if and only if for any two distinct indices $i,j$ in the same block $B_r$, we have that $d_i^{-1} \sltz d_i$ and $d_j^{-1} \sltz d_j$ are commensurable. This in turn is equivalent to $x_i/x_j\in\Q$. 

Hence, if $x_i/x_j\notin\Q$ for all $i\neq j$, then $\mathcal P$ has to be the singleton partition, and thus $\tilde H_{\mathcal P}=\sltr^n$. This proves the first implication.

Suppose on the other hand that $x_i/x_j\in\mathbb Q$ for some $i\ne j$, and let
\begin{equation}
\mathcal P_0=\{\{i,j\}\}\cup\{\{k\} : k\notin\{i,j\}\}.
\end{equation}
By the commensurability criterion above, $C=\Gamma\backslash\Gamma\psi(DH_{\mathcal P_0}D^{-1})$
is a proper closed orbit in $\Gamma\backslash\Gamma\psi(\mathrm{SL}(2,\mathbb R)^n)$.
Since $\varphi_X(\mathrm{SL}(2,\mathbb R))\subseteq\psi(DH_{\mathcal P_0}D^{-1})$, its orbit closure satisfies
\begin{equation}
\mathcal O_X\subseteq C
\subsetneq\Gamma\backslash\Gamma\psi(\mathrm{SL}(2,\mathbb R)^n).
\end{equation}
We conclude $m_X$ cannot be the invariant probability measure on the full product orbit. This proves the converse.
\end{proof}

We now obtain the limit theorem for the theta sums. 

\begin{prop}\label{prop:one-parameter-limit-law}
Let $f\in\Si(\R^n)$ and $X\in\Sym_n(\R)$ with $\det X\neq 0$.  For any $F \in C_b (\R_{\geq 0})$ and any absolutely continuous Borel probability measure $\lambda$ on $\R$, we have
    \begin{align}\label{eq:one-parameter-limit-law}
        \lim_{M \to \infty}\int_{\R} F\lp M^{-n/2}|S_f(M,tX)| \rp \de \lambda(t) = \int_{\oGam\bs \spdr} F(|\Theta_f(\oGam g)|) \de m_X(\oGam g) .
    \end{align}
\end{prop}

\begin{proof}
We apply Theorem \ref{cor:application-of-shah} with $\Gamma=\overline\Gamma_\Theta$, which is commensurable with $\Sp(n,\Z)$.
In particular, defining $\Psi_{f,X}(R)=m_X\{|\Theta_f|>R\}$ for every
$R\in\R$ gives $\Psi_{f,X}(R)=1$ for $R<0$. The Portmanteau theorem
then proves \eqref{LL003} at its continuity points.
\end{proof}

\section{Riemann-integrable cut-offs}
\label{sec:second}\label{sec:sharp}

We now extend the limit theorems to bounded, rapidly decreasing
functions that are Riemann-integrable on every compact rectangle.
Here rapid decrease means that for every $A>0$ there is $C_A<\infty$
such that $|f(\vecw)|\leq C_A(1+\|\vecw\|)^{-A}$.
This includes characteristic functions of bounded sets whose boundaries
have Lebesgue measure zero.

\begin{proposition}\label{propVar2}
Let $\lambda \ll m_\Theta$ be a Borel probability measure  with density $\rho\in L^\infty(\T_\Theta)$.  For every bounded, rapidly decreasing, Riemann-integrable function $f:\R^n\to\C$,
\begin{align}
\lim_{M\to\infty}\frac1{M^n}\int_{\T_\Theta}|S_f(M,X)|^2\,\d\lambda(X)
&=\|f\|_2^2+\int_{\R^n}f(\vecx)\overline{f(-\vecx)}\,\d\vecx \notag\\
&=2\|P_+f\|_2^2.
\label{eq:second-moment-full}
\end{align}
In particular, if $f$ is even, the limit equals $2\|f\|_2^2$.
\end{proposition}

\begin{proof}
For $f\in\Si(\R^n)$ this is the case $f_1=f_2=f$ of Theorem~\ref{thm:theta-sum-variance}. For general $f$, choose a sequence $f_j\in\Si(\R^n)$ with $\|f_j-f\|_2\to0$ as $j \to \infty$. Applying Proposition \ref{prop:covariance-limit} gives
\begin{equation}\label{eq:sharp-full-L2-approximation}
\limsup_{M\to\infty}\frac1{M^n}\int_{\T_\Theta}|S_f(M,X)-S_{f_j}(M,X)|^2\,\d\lambda(X)
\leq 2 \|\rho\|_\infty\|f-f_j\|_2^2.
\end{equation}
The result then follows.

\end{proof}

\begin{proposition}\label{prop007S}
Let $f:\R^n\to\C$ be bounded, rapidly decreasing, and
Riemann-integrable.  Let $\lambda\ll m_\Theta$ be a Borel probability measure on $\T_\Theta$. Then, for every $F\in C_b(\R_{\geq0})$,
\begin{equation}\label{eq:sharp-full-limit}
\lim_{M\to\infty}\int_{\mathbb T_\Theta} F\!\left(M^{-n/2}|S_f(M,X)|\right)\,\d\lambda(X)
=\int_{\overline\Gamma_\Theta\bs\spdr}F(|\Theta_f|)\,\d\mu_\Theta,
\end{equation}
where $\Theta_f$ denotes the $L^2$ extension given by Proposition \ref{L2}.
\end{proposition}

\begin{proof}
Choose $f_j\in\Si(\R^n)$ with $\|f-f_j\|_2\to0$. By the triangle inequality we have that
    \begin{align}
    \label{eq:14}
    &\Bigg|\int_{\T_\Theta} F\!\left(M^{-n/2}|S_f(M,X)|\right)\,\d\lambda(X) - \int_{\overline\Gamma_\Theta\bs\spdr}F(|\Theta_f|)\,\d\mu_\Theta\Bigg|\\
    &\leq  \int_{\T_\Theta} \Big|F\!\left(M^{-n/2}|S_{f}(M,X)|\right) - F\!\left(M^{-n/2}|S_{f_j}(M,X)|\right)\Big|\,\d\lambda(X)  \\
    & +\:\Bigg|\,\int_{\T_\Theta} F\!\left(M^{-n/2}|S_{f_j}(M,X)|\right)\,\d\lambda(X) - \int_{\overline\Gamma_\Theta\bs\spdr}F(|\Theta_{f_j}|)\,\d\mu_\Theta\Bigg|  \\
    & +\int_{\overline\Gamma_\Theta\bs\spdr}\Big|F(|\Theta_f|)- F(|\Theta_{f_j}|) \Big|\,\d\mu_\Theta.
  \end{align}
First suppose that $F$ is bounded and has Lipschitz constant $L$.   If $\lambda$ has density $\rho \in L^\infty$, then by Cauchy--Schwarz and \eqref{eq:sharp-full-L2-approximation} we have
\begin{multline}
  \limsup_{M\to\infty}\int_{\T_\Theta} \left|F(M^{-n/2}|S_f(M,X)|)-F(M^{-n/2}|S_{f_j}(M,X)|)\right|\,\d\lambda \\
  \leq L\sqrt{2\|\rho\|_\infty}\,\|f-f_j\|_2.
\end{multline}
Similarly, Proposition~\ref{L2} gives
\begin{equation}
\int_{\oGam\bs\spdr}|F(|\Theta_f|)-F(|\Theta_{f_j}|)|\,\d\mu_\Theta
\leq L\sqrt2\,\|f-f_j\|_2.
\end{equation}
For each fixed $j$, Theorem~\ref{thm:full-moment-convergence} shows that given $\veps > 0$, we may take $M$ sufficiently large so that
\begin{equation}
  \Bigg|\,\int_{\T_\Theta} F\!\left(M^{-n/2}|S_{f_j}(M,X)|\right)\,\d\lambda(X) - \int_{\overline\Gamma_\Theta\bs\spdr}F(|\Theta_{f_j}|)\,\d\mu_\Theta\Bigg| < \veps.
\end{equation}
Therefore, taking  $j$ sufficiently large and then $M$ sufficiently large gives the result for bounded Lipschitz test functions $F$.
 The extension to general $F \in C_b(\R_{\geq 0})$ follows from a standard argument.
 Furthermore, we may  approximate an arbitrary probability density $\rho$ in $L^1$ by bounded probability densities. Since $F$ is bounded, the resulting error is uniform in $M$, completing the proof.
\end{proof}

For diagonal forms recall the coordinatewise-even projection
$P_\square$ as in \eqref{eq:parity-projections}.

\begin{proposition}\label{propVar2DIAG}
Let $\lambda \ll m_\Theta^\Diag$ be a Borel probability measure on $\T_\Theta^{\Diag}$ with density $\rho\in L^\infty(\T_\Theta^\Diag)$.  Then, for every bounded, rapidly decreasing, Riemann-integrable $f:\R^n\to\C$,
\begin{align}
  \lim_{M\to\infty}\frac1{M^n}\int_{\T_\Theta^\Diag} |S_f(M,X)|^2\,\d\lambda(X)
&=\sum_{\boldsymbol\veps\in\{\pm1\}^n} \int_{\R^n}f(\vecw)\overline{f(\boldsymbol\veps\vecw)}\,\d\vecw \notag\\
&=2^n\|P_\square f\|_2^2.
\label{eq:second-moment-diagonal}
\end{align}
In particular, if $f$ is invariant under all coordinate reflections,
the limit equals $2^n\|f\|_2^2$.
\end{proposition}

\begin{proof}
If $f$ is Schwartz, then this is simply  Theorem~\ref{thm:diagonal-covariance} with
$f_1=f_2=f$. For Riemann-integrable $f$, we may approximate $f$ with Schwartz cut-offs $f_j$ such that
$\|f-f_j\|_2\to0$. By \eqref{eq:diag-covariance-limit} we have
\begin{equation}
\limsup_{M\to\infty}M^{-n}\int_{\T_\Theta^\Diag}  |S_f(M,X)-S_{f_j}(M,X)|^2\,\d\lambda \leq 2^n\|\rho\|_\infty\|f-f_j\|_2^2.
\end{equation}
The triangle inequality in $L^2$ and continuity of $P_\square$ then
give \eqref{eq:second-moment-diagonal}.
\end{proof}

\begin{proposition}
Let $f:\R^n\to\C$ be bounded, rapidly decreasing, and Riemann-integrable.  Let $\lambda\ll m_\Theta^\Diag$ be a Borel probability measure on $\T_\Theta^\Diag$.  Then, for every $F\in C_b(\R_{\geq0})$,
\begin{align}
\lim_{M\to\infty}\int_{\T_\Theta^\Diag} F\!\left(M^{-n/2}|S_f(M,X)|\right)\, \d\lambda(X) =\int_{\oDel\bs(\sltr)^n} F(|\Theta_f|)\, \d\dlimmeas.
\label{eq:sharp-diagonal-limit}
\end{align}
\end{proposition}

\begin{proof}
Arguing as  in the proof of  Proposition~\ref{prop007S}, with  Proposition~\ref{propVar2DIAG} and \eqref{eq:L2-product} in place of  Proposition~\ref{prop:covariance-limit} and Proposition~\ref{L2} gives the result. 
\end{proof}

\section{Tail asymptotic for generic \texorpdfstring{$X$}{X}}
\label{section:Rank-n-asymptotics}

For every $R\in\R$, define
\begin{equation}\label{eq:generic-tail-definition}
\Psi_f(R):=\mu_\Theta\{\oGam g:|\Theta_f(g)|>R\};
\end{equation}
in particular $\Psi_f(R)=1$ for $R<0$.
This is a right-continuous, non-increasing function. Theorem
\ref{thm:full-moment-convergence} and the Portmanteau theorem give
\eqref{LL001} at every continuity point, equivalently at every $R$
with $\mu_\Theta\{|\Theta_f|=R\}=0$.
There are at most countably many exceptional values.

We  define
\begin{equation}\label{eq:volume-products-intro}
\mathcal V_n:=\prod_{j=1}^n\frac{\Gamma(j)\zeta(2j)}{\pi^j},
\qquad
\mathcal P_n:=\prod_{j=2}^n
\pi^{-j/2}\Gamma\!\lp\frac j2\rp\zeta(j).
\end{equation}
Let $\omega_n$ be Haar probability measure on $\Un$, and let
\begin{equation}\label{eq:generic-angular-moment}
\mathcal M_{n,f}:=\int_{\Un}|f_Q(\vecnull)|^{2(n+1)}\,\d\omega_n(Q).
\end{equation}

\begin{theorem}\label{thm:generic-tail}
Let $n\geq2$ be fixed and let $f\in\Si(\R^n)$.  For every
$\varepsilon>0$,
\begin{equation}\label{eq:generic-tail}
\Psi_f(R)=C_{n,f}R^{-2(n+1)}
+O_{n,f,\varepsilon}\!\left(R^{-2(n+2)+\varepsilon}\right),
\end{equation}
where
\begin{equation}\label{eq:generic-tail-constant}
C_{n,f}=\frac{2\mathcal P_n}
{(n+1)(2^n+1)\mathcal V_n}\,\mathcal M_{n,f}.
\end{equation}
The coefficient is positive if and only if $P_+f\neq0$.
\end{theorem}

We recall that $\spdr$ is naturally identified with $\uhp_n \times \Un$ via the Iwasawa decomposition. Up to sets of measure zero, a fundamental domain for $\oGam \bs \spdr$ is identified with $\calF_\Theta \times K_0 \subseteq \uhp_n \times \Un$. The set $\calF_\Theta$ is in fact a  fundamental domain for $\oGam$ on $\uhp_n$, where $\oGam$ acts as in \eqref{def:siegel-mobius}, and $K_0$ is a fundamental domain for $\{\pm I\}\bs\Un$.

If $\calF_n$ is a fundamental domain for the action of $\spdz$ on $\uhp_n$ then we may choose $\calF_\Theta$ to be 
\begin{equation}
  \label{eq:fund-dom-gamma-theta}
  \calF_\Theta := \bigcup_{\gamma \in \Lambda_\Theta}\gamma\calF_n 
\end{equation}
where $\Lambda_\Theta$ denotes a complete set of coset representatives for $\overline{\Gamma}_\Theta \bs \spdz$. Let $\mu$ denote the Haar measure on $\spdr$. In `Iwasawa coordinates' the  measure $\mu$ is given by
\begin{equation}
  \label{eq:haar-iwasawa-coordinates}
 \de \mu =  \frac{\de X \de Y}{(\det Y)^{n+1}} \de \omega_n(Q)
\end{equation}
where $\de X$ denotes Lebesgue measure on $\Sym_n(\R)$,  $\de Y$ denotes Lebesgue measure on the space of symmetric, positive definite matrices $\scrP_n$, and $\omega_n$ is Haar probability measure on $\Un$. If we write $Y = UV \trans U$, with $U$ unipotent, upper triangular, and $V = \diag (\vecv)$, $v_i > 0$, then this measure becomes
\begin{equation}\label{eq:haar-iwasawa-measure}
     \de X  \de U \lp\prod_{1 \leq  j \leq n} v_j^{-n + j - 2} \de v_j\rp  \d\omega_n(Q),
\end{equation}
where $\de U$ and $\de v_j$ are the Lebesgue measures on the entries of $U$, and $\boldsymbol v$, respectively.
 The corresponding volume element on $\uhp_n$, invariant under the action \eqref{def:siegel-mobius} is given by
\begin{equation}\label{eq:siegel-volume-element}
\d\vol_{\uhp_n}(Z) = \frac{\d X\,\d Y}{(\det Y)^{n+1}},
\end{equation}
so that $\d\mu(Z,Q)=\d\vol_{\uhp_n}(Z)\,\d\omega_n(Q)$. A classical result due to Siegel \cite{Siegel-Symplectic} is that 
\begin{equation}\label{volFD0}
\vol_{\uhp_n}(\calF_n) = 2 \prod_{j=1}^n \frac{\Gamma(j) \zeta(2j)}{\pi^j},
\end{equation}
and so
\begin{equation}\label{volFD}
\mathcal{V}_n = \mu(\calF_{n} \times K_0) =  \prod_{j=1}^n \frac{\Gamma(j) \zeta(2j)}{\pi^j},
\end{equation}
where the factor $\frac12$ is because  $\omega_n(K_0) = \tha$.
Through the group law of $\Gamma_\Theta$ given in \eqref{eq:Gamma-explicit-definition} we define an affine action on $(\Z / 2\Z)^{2n}$  by
\begin{equation}
  \label{eq:Gamma-Theta-action-on-torus}
  \vecv \mapsto \vecv\gamma\inv + 2\vecp_\gamma
\end{equation}
where $\vecp_\gamma = (\vecr, \vecs) \in \R^n \times \R^n$, with 
\begin{equation}
  2r_i\equiv(C\trans D)_{ii} \pmod2, \qquad 2s_i \equiv(A\trans B)_{ii}\pmod2.
\end{equation}

The orbit $\calE_n$ of $(\bn,\bn)$ is the set of
\emph{even theta characteristics}:
\begin{equation}
\calE_n=\{(\boldsymbol a,\boldsymbol b)\in(\Z/2\Z)^n\times(\Z/2\Z)^n:
\boldsymbol a\boldsymbol b^{\mathrm t}=0\text{ in }\Z/2\Z\}.
\end{equation}
The bijection is given by the map $\oGam\gamma\mapsto2\vecp_{\gamma^{-1}}$. The stabiliser of $(\bn, \bn)$ is $\oGam$. Therefore $|\calE_n| = [\spdz : \oGam] = 2^{n-1}(2^n + 1)$, see \cite[pp. 226, 235]{Igusa-theta-constants} for details. It then follows that
\begin{equation}
  \mu(\calF_\Theta \times K_0) = |\calE_n| \mathcal{V}_n = 2^{n - 1}(2^n + 1) \prod_{j=1}^n \frac{\Gamma(j)\zeta(2j)}{\pi^j}.
\end{equation}
Consequently, we have that
\begin{equation}
\Psi_f(R)=\frac1{|\mathcal E_n|\mathcal V_n}
\sum_{\boldsymbol m\in\mathcal E_n}I_{\boldsymbol m}(R),
\quad
I_{\boldsymbol m}(R):=\int_{\calF_{n}\times K_0}
\ind_{\{|\Theta_f(\boldsymbol m/2;g)|>R\}}\,\d\mu(g).
\end{equation}

Let
\begin{equation}
\calE_n^0 := \{\bn\}\times(\Z/2\Z)^n\subset\calE_n,\qquad |\calE_n^0|=2^n.
\end{equation}

We then further decompose $\Psi_f$ as  
\begin{equation}
  \Psi_f(R) = \frac{1}{|\calE_n| \mathcal V_n} \lp\sum_{\vecm \in \calE^{0}_n} I_{\vecm}(R) + \sum_{\vecm \in\calE_n \bs  \calE^{0}_n} I_{\vecm}(R)\rp.
\end{equation}

A convenient choice of fundamental domain for $\spdz$ was  constructed in \cite{Marklof-Welsh-higher-rank-theta-sumsII}. We give the relevant definitions below.

Let
\begin{equation}
\calF_1:=\{x+\i y\in\uhp_1:|x|\leq\tfrac12,\ x^2+y^2\geq1\}.
\end{equation}

For $n\geq 2$ and $g\in \Sp(n,\R)$ write
\begin{equation}
  g=
  \begin{pmatrix}
    1&\boldsymbol{r}_1&t_1-\boldsymbol{s}_1 \trans\boldsymbol{r}_1&\boldsymbol{s}_1\\
    0&I&\trans\boldsymbol{s}_1&0\\
    0&0&1&0\\
    0&0&-\trans\boldsymbol{r}_1&I
  \end{pmatrix}
  \begin{pmatrix}
    1&0&0&0\\ 0&I&0&X_1\\ 0&0&1&0\\ 0&0&0&I
  \end{pmatrix}
  \begin{pmatrix}
    v_1^{1/2}&&&\\ &Y_1^{1/2}&&\\ &&v_1^{-1/2}&\\ &&&\trans Y_1^{-1/2}
  \end{pmatrix}
  k(Q),
\end{equation}
where $\boldsymbol{r}_1,\boldsymbol{s}_1\in\R^{n-1}$, $t_1\in\R$,
$X_1$ symmetric, $v_1>0$, $Y_1$ positive definite symmetric, $Q\in\Un$,
and set
\begin{equation}
  g_1(g)=\begin{pmatrix}I&X_1\\0&I\end{pmatrix}
         \begin{pmatrix}Y_1^{1/2}&0\\0&\trans Y_1^{-1/2}\end{pmatrix}
         \in\Sp(n-1,\R).
\end{equation}
Then $\calF_n\subset\uhp_n$ consists of the points $Z=X+\i Y$ whose representative $g=n(X)a(Y)$ satisfies
\begin{enumerate}
  \item $v_1(g)\ge v_1(\gamma g)$ for all $\gamma\in\spdz$;
  \item $g_1(g)\cdot(\i I_{n-1})\in\calF_{n-1}$;
  \item $|(\boldsymbol{r}_1)_j|,\,|(\boldsymbol{s}_1)_j|,\,|t_1|\le\tfrac12$
        for all indices, with $(\boldsymbol{r}_1)_1\ge0$.
\end{enumerate}

We recall that given $Y$ positive definite, we have the decomposition $Y = U V \trans U$, where $U$ is an upper triangular unipotent matrix, and $V = \diag(v_1, \dots, v_n)$. Given $X + \i Y \in \calF_n$ we have by \cite[Proposition 3.2]{Marklof-Welsh-higher-rank-theta-sumsII} that for $1 \leq j \leq n -1$
\begin{equation}
  \label{eq:fd-reduction}
  \qquad v_j\geq\frac34v_{j+1},\quad v_n\geq \frac{\sqrt 3}{2}.
\end{equation}

Central to proof of Theorem~\ref{thm:generic-tail} are the estimates \eqref{eq:global-theta-bound-tail}--\eqref{eq:general-theta-asymptotic}, which are specific instances of the more general expansion given in \cite[Theorem 4.4]{Marklof-Welsh-higher-rank-theta-sums}. Typically such estimates are stated within a fixed fundamental domain, but in fact hold in a larger region which is much easier to integrate over. To this end, fix $B_n<\infty$ large enough to bound the entries of $X$ and the
off-diagonal entries of $U$ in $\calF_n$, and in $\calS_P$ defined in \eqref{eq:parabolic-section} below. For $b\geq0$ define
\begin{equation}\label{eq:bounded-reduction-region}
\begin{split}
\mathcal R_n(b):=\{X+\i UV\trans U:\ &|X_{ij}|,|U_{ij}|\leq B_n \text{ for }
i \leq j\\
&U\text{ upper unipotent},\
v_j\geq\tfrac34v_{j+1} \text{ for }j<n,\ v_n>b\}.
\end{split}
\end{equation}
Here and below $V$ is a diagonal matrix with positive entries. We have for any $\vecw \in \R^n$, uniformly on $\mathcal R_n(0)$,
\begin{equation}\label{eq:reduction-quadratic-comparison}
\vecw Y\trans\vecw\asymp_n\vecw V\trans\vecw
\end{equation}
Indeed, $V^{-1/2}UV^{1/2}$ and its inverse have bounded entries.
For fixed $b>0$, \eqref{eq:theta-abs-conv} and the uniform Schwartz
estimates therefore imply, for $Z\in\mathcal R_n(b)$,
$Q\in K_0$, $\vecx\in[-1/2,1/2]^n$, $\vecy\in\R^n$, and every $A>0$,
\begin{equation}\label{eq:global-theta-bound-tail}
|\Theta_f((\vecx,\vecy),Z,Q)|
\ll_{n,f,A,b}(\det V)^{1/4}(1+\vecx V\trans\vecx)^{-A}.
\end{equation}
The bound is uniform since $(m_j+x_j)^2\gg m_j^2+x_j^2$ for $|x_j|\leq1/2$, and so 
we may take $A$ sufficiently large before summing over $\vecm$.

For $\vecm\in\calE_n^0$, we have uniformly in $Z\in\mathcal R_n(0)$ and $Q\in K_0$,
\begin{equation}\label{eq:39}
\left||\Theta_f(\vecm/2,Z,Q)|-(\det Y)^{1/4}|f_Q(\bn)|\right|
\ll_{n,f,L}(\det Y)^{1/4}
\sum_{\vecb\in\Z^n\setminus\{\bn\}}(1+\vecb V\trans\vecb)^{-L}.
\end{equation}
If $v_n\gg_n T\geq1$, then
$\vecb V\trans\vecb\gg_n T\|\vecb\|^2$. Given $A>0$, we set $L>\max\{A,n/2\}$ so we have the estimate
\begin{equation}\label{eq:general-theta-asymptotic}
\left||\Theta_f(\vecm/2,Z,Q)|-(\det Y)^{1/4}|f_Q(\bn)|\right|
\ll_{n,f,A}(\det Y)^{1/4}T^{-A},\qquad \vecm\in\calE_n^0.
\end{equation}

In particular, these estimates suggest that $|\Theta_f|$ can only be large in a certain region of the cusp of $\calF_n$, which we now describe precisely. Let $\scrP_n$ denote the positive definite symmetric $n\times n$ matrices.
Set $\calD_1=\{y>0\}$ and

\begin{equation}
  \calD_2=\left\{
  \begin{pmatrix}1&r\\0&1\end{pmatrix}
  \begin{pmatrix}v_1&0\\0&v_2\end{pmatrix}
  \begin{pmatrix}1&0\\r&1\end{pmatrix}
  \;:\;v_1,v_2>0,\ 0\le r\le\tfrac12,\ r^2+\frac{v_1}{v_2}\ge1 \right\}.
\end{equation}

For $l>2$, writing $Y\in\scrP_l$ as
\begin{equation}
  Y=\begin{pmatrix}1&\boldsymbol{r}\\0&I\end{pmatrix}
    \begin{pmatrix}v_1&0\\0&Y_1\end{pmatrix}
    \begin{pmatrix}1&0\\ \trans\boldsymbol{r}&I\end{pmatrix},
  \qquad \boldsymbol{r}\in\R^{l-1},\ v_1>0,\ Y_1\in\scrP_{l-1},
\end{equation}
define $\calD_l$ inductively as the set of such $Y$ with
\begin{enumerate}
  \item $v_1(Y)\ge v_1\!\left(A\,Y\,\trans A\right)$ for all $A\in\GL(l,\Z)$;
  \item $Y_1\in\calD_{l-1}$;
  \item $|r_j|\le\tfrac12$ for all $j$, and $0\le r_1\le\tfrac12$.
\end{enumerate}
Equivalently, condition (1) states that
\begin{equation}
  \frac{1}{v_1(Y)}=\min\bigl\{\,\boldsymbol{c}\,Y^{-1}\,\trans\boldsymbol{c}
  \;:\;\boldsymbol{c}\in\Z^{l}\setminus\{0\}\,\bigr\},
\end{equation}
since $v_1(Y)^{-1}=(Y^{-1})_{11}$.

Let $P<\spdr$ be the parabolic
\begin{equation}
  \label{eq:Siegel-parabolic}
  P := \left\{\matr{A}{S\trans A^{-1}}{0}{\trans A^{-1}} : A\in\GL(n,\R),\ S\in\Sym_n(\R) \right\}
\end{equation}
and $\Gamma_P=P\cap\spdz$. 
 We define the invariant parabolic height as
\begin{equation}\label{eq:parabolic-height}
h_P(\Gamma_Pg):=\min_{A\in\mathrm{GL}(n,\Z)}v_n(AY\trans A).
\end{equation}
By \cite[Lemma 3.3(3)]{Marklof-Welsh-higher-rank-theta-sumsII},
\begin{equation}\label{eq:height-comparison}
h_P(\Gamma_Pg)\asymp_n v_n(Y),
\qquad Y\in\calD_n.
\end{equation}
Let $\mathcal X_n=[-1/2,1/2)^{n(n+1)/2}$.  Up to a set of $\mu$-measure zero, the set $\mathcal S_P \times K_0$, where
\begin{equation}\label{eq:parabolic-section}
\mathcal S_P=\{X+\i Y:Y\in\calD_n,\ X\in\mathcal X_n\},
\end{equation}
is a fundamental domain for $\Gamma_P\bs\spdr$, under the
identification $(X+\i Y,Q)\leftrightarrow n(X)a(Y)k(Q)$.
Both $\mathcal S_P$ and $\calF_n$ lie in $\mathcal R_n(0)$, and
$\calF_n\subset\mathcal R_n(b_0)$ for any fixed
$0<b_0<\sqrt3/2$. On $\mathcal S_P$, however, $v_n$ can tend to
zero. By \eqref{eq:height-comparison}, its region $h_P>T$ lies in
$\mathcal R_n(b_0)$ once $T$ is sufficiently large.
The comparison \eqref{eq:height-comparison} also holds on $\calF_n$,
whose imaginary part belongs to $\calD_n$. An element of $\Gamma_P$ as in \eqref{eq:Siegel-parabolic} acts on $\uhp_n$ via \eqref{def:siegel-mobius}, that is $Z\mapsto AZ\trans A+S$, with $A\in\mathrm{GL}(n,\Z)$ and $S\in\Sym_n(\Z)$.  Using Grenier reduction, followed by inversion, we send  $Y$ into $\calD_n$ and if the image is not on the boundary, then $A$ is unique up to sign. Translation maps $X$ into $\mathcal X_n$.

\begin{lemma}\label{lem:deep-cusp}
There is $H_0(n)<\infty$ such that the  map
\begin{align}\label{eq:deep-cusp-map}
  \Gamma_P\bs\spdr &\longrightarrow\spdz\bs\spdr,\\
  \Gamma_P g &\longmapsto \spdz g.
\end{align}
is injective on $\{h_P>H_0(n)\}$.  Furthermore, for $H_0$ sufficiently large, we have that for every $H\geq H_0$ and every nonnegative or integrable $\Gamma_P$-invariant function $F$,
\begin{equation}\label{eq:deep-cusp-unfolding}
\int_{\mathcal{F}_n\times K_0}F(g)\ind_{\{h_P>H\}}\,\d\mu(g)
=\int_{\mathcal S_P\times K_0}F(g)\ind_{\{h_P>H\}}\,\d\mu(g).
\end{equation}

\end{lemma}

\begin{proof}
Note that if $\veca$ is the last row of $A$, then $v_n(AY\trans A)=\veca Y\trans \veca$, and as $A \in \mathrm{GL}(n, \Z)$, the possible last rows are exactly the primitive integral vectors, $\Z^n_{\mathrm{prim}}$.  Hence
\begin{equation}\label{eq:primitive-height}
h_P(\Gamma_Pg)=
\min_{\veca\in\Z^n_{\mathrm{prim}}}\veca Y\trans \veca,
\end{equation}
and so $\vecc Y\trans \vecc\geq h_P$ for every nonzero $\vecc\in\Z^n$.
For reduced $Y\in\calD_n$, the same estimates give, for some $\kappa_n>0$,
\begin{equation}\label{eq:eigenvalue-height}
\lambda_{\min}(Y)\geq\kappa_n h_P(\Gamma_Pg).
\end{equation}
  
Consider points $\Gamma_P g_1, \Gamma_P g_2$, for which $\spdz g_1 = \spdz g_2$. Let $h_j := h_P (\Gamma_P g_j)$. Choose representatives $Z_j = X_j + \i Y_j \in \uhp_n$ for $g_j$ with $Y_j\in\calD_n$ and $Z_2=\gamma Z_1$, where $\gamma=\smatr{A}{B}{C}{D}\in\spdz$.  The imaginary part transforms as 
\begin{equation}\label{eq:imaginary-inverse}
Y_2^{-1}=(CX_1+D)Y_1^{-1}\trans (CX_1+D)+CY_1\trans C.
\end{equation}
If $C\ne0$, choose a nonzero integral row $\vecc$ of $C$.  Since $Y_1\inv$ is positive definite, we have that $(CX_1+D)Y_1^{-1}\trans (CX_1+D)$ is positive semidefinite, and so
\begin{equation}\label{eq:height-lower-bound-intermediate}
  \lambda_{\min}(Y_2)^{-1} = \lambda_{\max} (Y_2^{-1}) \geq \vecc Y_1\trans \vecc\geq h_1.
\end{equation}
By \eqref{eq:eigenvalue-height}, it follows that
\begin{equation}\label{eq:height-separation}
h_1h_2\leq\kappa_n^{-1}.
\end{equation}
Thus $C=0$ if both $h_1$ and $h_2$ are greater than $\kappa_n^{-1/2}$. As $\gamma \in \spdz$, we have $A\in\mathrm{GL}(n,\Z)$, $D=\trans A^{-1}$, and $B\trans A$ symmetric. Therefore $\gamma\in\Gamma_P$, and the map is injective.

Finally, by  \cite[Proposition 3.5, with $l=n$]{Marklof-Welsh-higher-rank-theta-sumsII} we may take $H_0$ sufficiently large so that, whenever $h_P(\Gamma_P g)\geq H_0$, there exists $\gamma \in \Gamma_P$ such that $\gamma g \in \calF_n \times K_0$, and  hence \eqref{eq:deep-cusp-unfolding}.
\end{proof}

The following lemma is a general volume estimate that will be useful in the proofs of Lemmata \ref{lem:tail-sum-reduction} and \ref{lem:low-cusp-bound}; cf.~the proof of Lemma 5.1 in \cite{Marklof-Welsh-higher-rank-theta-sums}.

\begin{lemma}\label{lem:ordered-cone}
Let $1\leq m\leq n$, $0<c\leq1$, and $b,H>0$. Let $\lambda_j=n-j+1$, and define
\begin{equation}\label{eq:beta-m}
\beta_m:=\frac1m\sum_{j=1}^m\lambda_j
=n+1-\frac{m+1}{2}.
\end{equation}
Then
\begin{equation}\label{eq:ordered-cone}
\begin{aligned}
&\int_{\substack{v_m\geq b,\ v_j\geq cv_{j+1}\ (j<m)\\
v_1\cdots v_m\geq H}}
\prod_{j=1}^m v_j^{-\lambda_j-1}\,\d v_1 \cdots \de v_m\\
&\hspace{2cm}\ll_{n,m,c}\min\{b^{-m\beta_m},H^{-\beta_m}\}.
\end{aligned}
\end{equation}
\end{lemma}

\begin{proof}
Set $v_j=b\,c^{m-j}w_j$. Then $w_1\geq\cdots\geq w_m\geq1$,
the measure acquires a factor $C_{n,m,c}b^{-m\beta_m}$, and the
product threshold becomes $H'=H/(b^mc^{m(m-1)/2})$.
It therefore suffices to prove the normalised estimate
$O_{n,m}(\min\{1,(H')^{-\beta_m}\})$; the powers of $b$ cancel
in the second term. Relabel $w_j$ as $v_j$ and $H'$ as $H$.  Let $v_j=\e^{u_j}$ and set $r_i=u_i-u_{i+1}$ for $i<m$, and $r_m=u_m$. We have that $r_i\geq0$, and $\log(v_1\cdots v_m)=\sum_{i=1}^m i r_i$. Therefore
\begin{equation}\label{eq:ordered-cone-log-coordinates}
\prod_{j=1}^m v_j^{-\lambda_j-1}\,\de v_1 \cdots \de v_m =\exp\!\left(-\sum_{i=1}^m A_i r_i\right)\de r_1 \cdots \de r_m,
\end{equation}
where $A_i=\sum_{j=1}^i\lambda_j$.  A direct calculation gives
\begin{equation}\label{eq:ordered-cone-gap}
A_i-i\beta_m=\frac{i(m-i)}2>0\quad(i<m),
\qquad A_m=m\beta_m.
\end{equation}
If $L=\log H\geq0$, integration first in $r_m$ under
$\sum ir_i\geq L$ leaves exponential decay in the other variables and
gives $O(\e^{-\beta_mL})$.  If $H<1$, the integral is then over the entire region $r_1, \dots, r_m \geq 0$ which is finite.
\end{proof}

\begin{lemma}\label{lem:tail-sum-reduction}
Fix $b>0$, $\vecx\neq\bn$, and $A>n/4$. For $R\geq 1$,
\begin{equation}\label{eq:7}
\mu\{(Z,Q)\in\mathcal R_n(b)\times K_0:
(\det V)^{1/4}(1+\vecx V\trans\vecx)^{-A}>R\}
\ll_{n,\vecx,A,b}R^{-2(n+2)} .
\end{equation}
\end{lemma}

\begin{proof}
Let $k$ be the first index with $x_k\neq0$. On $\mathcal R_n(b)$,
we have $v_j\ll_n v_k$ for $j\geq k$, and all $v_j$ have fixed
positive lower bounds. If $k=1$, the expression in \eqref{eq:7}
is bounded because $A>n/4$, so its tail is eventually empty.
If $k>1$, the tail condition implies
\begin{equation}
v_1\cdots v_{k-1}\gg_{n,\vecx,A}R^4
v_k^{4A-(n-k+1)}\gg_{n,\vecx,A,b}R^4.
\end{equation}
By the decomposition of Haar measure given in \eqref{eq:haar-iwasawa-measure}, we may integrate out each of the $X,U,Q$ coordinates, as they are constrained to compact regions, and bound each integral over $v_j$, with $j\geq k$. The measure is bounded, up to a constant, by
\begin{equation}
  \int_{\substack{v_{k-1}\geq (3/4)^{n-k+1}b,\ v_j\geq cv_{j+1}\ (j<k-1)\\
v_1\cdots v_{k-1}\gg_{n, \vecx, A, b} R^4}}
\prod_{j=1}^{k-1} v_j^{-(n - j +2)}\,\d v_1 \cdots \de v_{k-1},
\end{equation}
and so Lemma~\ref{lem:ordered-cone} with $c = \tfrac34$, and $m=k-1$ bounds what remains by $O_{n,\vecx,A,b}(R^{-4(n+1)+2k})$. 
\end{proof}

The next lemma treats the part of the cusp where $v_n$ is bounded. 

\begin{lemma}\label{lem:low-cusp-bound}
For $T>0$ and $R>0$,
\begin{equation}\label{eq:low-cusp-bound}
\mu\{(Z,Q)\in\mathcal S_P\times K_0:
(\det V)^{1/4}>R,\ 0<v_n\leq T\}
\ll_n R^{-2(n+2)}T^{n/2}.
\end{equation}
The same bound holds with $\mathcal S_P$ replaced by
$\mathcal R_n(0)$, and hence by $\calF_n$.
\end{lemma}

\begin{proof}
It suffices to integrate over $\mathcal R_n(0)$. For fixed $v_n>0$, we apply Lemma~\ref{lem:ordered-cone} with $m=n-1$, $b=3v_n/4$, $c=3/4$ and $H=R^4/v_n$ to obtain \begin{equation}\label{eq:low-cusp-integral}
\begin{aligned}
&\mu\{(Z,Q)\in\mathcal R_n(0)\times K_0:
(\det V)^{1/4}>R,\ v_n\leq T\}\\
&\qquad\ll_n R^{-2(n+2)}\int_0^T v_n^{(n-2)/2}\,\d v_n
\ll_n R^{-2(n+2)}T^{n/2}.
\end{aligned}
\end{equation}
\end{proof}

Finally, before proceeding to the proof of Theorem \ref{thm:generic-tail}, we calculate the volume of the cusp of $\mathcal S_P$.

\begin{lemma}\label{lem:cusp-volume}
For $R>0$, we have
  \begin{equation}
    \label{eq:cusp-volume}
    \mu\{(X + \i Y, Q) \in \mathcal S_P \times K_0 :  \det Y > R\} = \frac{\mathcal P_n}{(n + 1)} R^{-(n+1)/2} .
  \end{equation}

\end{lemma}

\begin{proof}
  Let  $\calF_n^{\Gr}$ denote  Grenier's fundamental domain for the action of $\mathrm{GL}(n ,\Z)$ on $\scrP_n$. The \emph{Euclidean} volume of the following set can be shown to be (see for example, Grenier~\cite[Theorem I(3)]{Grenier})  
  \begin{equation}
    \label{eq:grenier-volume}
    \Leb\{ Y \in \calF_n^{\Gr} : \det Y < 1\} = \frac{2}{n + 1} \mathcal P_n,
  \end{equation}
  where $\mathcal P_n$ is as in \eqref{eq:volume-products-intro}.
  
  Let $\scrSP_n$ denote the determinant one surface in $\scrP_n$. Note that given $ Y \in \scrP_n$, then there exists $W \in \scrSP_n$ such that $Y = \rho W$, where $\rho =( \det Y)^{1/n}$.  The invariant measure on $\scrP_n$ disintegrates as follows, with $\d W$ normalised by this identity (cf.~\cite{Terras-VolII}): 
  \begin{equation}
    \label{eq:radial-disintegration}
    \frac{\de Y}{(\det Y)^{(n+1)/2}} = \rho^{-1} \,\de \rho \,\de W.
  \end{equation}
  It follows that
  \begin{equation}
    \label{eq:grenier-radial-integral}
   \Leb\{ Y \in \calF_n^{\Gr} : \det Y < 1\}  = \int_0^1 \rho^{\frac{n(n+ 1)}{2}-1} \, \de \rho \int_{\calF_n^{\Gr} \cap \scrSP_n} \, \de W,
 \end{equation}
 and so the measure of the determinant-one surface of the Grenier domain is
 \begin{equation}
   \label{eq:shape-volume}
   \int_{\calF_n^{\Gr} \cap \scrSP_n} \, \de W = n \mathcal P_n.
 \end{equation}

  The invariant measure is preserved by inversion, so the determinant-one surface of the domain $\calD_n$ also has measure $n\mathcal P_n$. Therefore
 \begin{align}
   & \label{eq:cusp-radial-integral}
   \mu\{(X + \i Y, Q) \in \mathcal S_P \times K_0 :  \det Y > R\} \notag \\
   & = \ha  n \mathcal P_n \int_{R^{\frac{1}{n}}}^\infty \rho^{-\tfrac{n(n+1)}{2} - 1} \, \de \rho = \frac{\mathcal P_n}{n+1} R^{-(n+1)/2}.
 \end{align}
\end{proof}

\begin{proof}[Proof of Theorem \ref{thm:generic-tail}]
We first note that
\begin{multline}
  (|\mathcal E_n|\mathcal V_n)\Psi_f(R) =\sum_{\vecm \in \calE_n}\mu\{g \in \calF_n \times K_0 : |\Theta_f (\vecm/2, g)| > R, h_P(g)  > R^{2\veps/n}\} + \\
  \sum_{\vecm \in \calE_n}\mu\{g \in \calF_n \times K_0 : |\Theta_f (\vecm/2, g)| > R, h_P( g)  \leq R^{2\veps/n}\}.
\end{multline}
Using \eqref{eq:global-theta-bound-tail}, followed by Lemma \ref{lem:low-cusp-bound}, with a constant multiple of $T=R^{2\veps/n}$, we have, since $h_P\asymp_n v_n$ on $\calF_n$,
  \begin{equation}
    \mu \{g \in \calF_n \times K_0 : |\Theta_f (\vecm/2, g)| > R, h_P( g)  \leq R^{2 \veps / n}\} \ll_{n,f} R^{-2(n + 2)+ \veps}.
  \end{equation}
  Using that the function
  \begin{equation}
    g \mapsto \sum_{\vecm \in \calE_n} \ind\{|\Theta_f (\vecm/2, g)| > R\}
  \end{equation}
  is $\spdz$-invariant, we have by Lemma \ref{lem:deep-cusp} that
  \begin{align}
    &\sum_{\vecm \in \calE_n}\mu\{g \in \calF_n \times K_0 : |\Theta_f (\vecm/2, g )| > R,\: h_P( g)  > R^{2\veps/n}\}\\
    &=\sum_{\vecm \in \calE_n}\mu\{g \in \mathcal{S}_P \times K_0 : |\Theta_f (\vecm/2, g )| > R,\: h_P( g)  > R^{2\veps/n}\}.
  \end{align}
  For $\vecm\in\calE_n\setminus\calE_n^0$, Lemma~\ref{lem:tail-sum-reduction} gives
  \begin{equation}
    \sum_{\vecm \not\in \calE^0_n}\mu\{g \in \mathcal{S}_P \times K_0 : |\Theta_f (\vecm/2, g )| > R,\: h_P( g)  > R^{2\veps/n}\} \ll_{n,f} R^{-2(n + 2)},
  \end{equation}
  so we are left to estimate
  \begin{equation}
    \sum_{\vecm \in \calE_n^0}\mu\{g \in \mathcal{S}_P \times K_0 : |\Theta_f (\vecm/2, g )| > R, h_P(g)  > R^{2\veps/n}\}.
  \end{equation}

  Let $p=2(n+1)$ and $T=R^{2\varepsilon/n}$, and fix $\vecm\in\calE_n^0$. Define 
\begin{equation}
a_-(Q)=\max\{|f_Q(\bn)|-C_{n,f,A}T^{-A},0\},\qquad a_+(Q)=|f_Q(\bn)|+C_{n,f,A}T^{-A}.
\end{equation}
Then we have 
 \begin{multline}
    \label{eq:34}
    \mu\{(Z,Q) \in \mathcal{S}_P \times K_0 : (\det Y)^{1/4}a_-(Q)> R,\: h_P(Z,Q)  > T\} \leq \\\mu\{(Z,Q) \in \mathcal{S}_P \times K_0 : |\Theta_f (\vecm/2, Z, Q )| > R,\: h_P(Z,Q)  > T\}\\ \leq \mu\{(Z,Q) \in \mathcal{S}_P \times K_0 : (\det Y)^{1/4}a_+(Q)> R,\: h_P(Z,Q)  > T\}.
  \end{multline}

The functions $a_\pm$ are bounded uniformly for $T\geq1$. We observe that removing the condition $h_P>T$ from the upper and lower bounds gives an error of at most $O_{n,f,A}(R^{-2(n+2)}T^{n/2})$, by Lemma~\ref{lem:low-cusp-bound}, since $\{h_P\leq T\}\subseteq\{v_n\leq c_nT\}$.
Applying Lemma~\ref{lem:cusp-volume} then gives, for both bounds,
\begin{align}
&\mu\{(Z,Q)\in\mathcal S_P\times K_0:(\det Y)^{1/4}a_\pm(Q)>R,\ h_P(Z,Q) >T\} \notag \\
&\quad=\frac{2\mathcal P_n}{n+1}R^{-p}
\int_{K_0}a_\pm(Q)^p\,\d\omega_n(Q)
+O_{n,f,\varepsilon}(R^{-p-2+\varepsilon}).
\end{align}
Since $|f_Q (\bn)|$ is bounded, $\int_{K_0}|a_\pm(Q)^p-|f_Q(\bn)|^p|\,\d\omega_n=O_{n,f,A}(T^{-A})$.
Taking $A$ so that $2A\varepsilon/n\geq2$ we have that
\begin{align}\label{eq:36}
&\mu\{g\in\mathcal S_P\times K_0:|\Theta_f(\vecm/2,g)|>R,\ h_P(g)>T\}
\notag\\
& =\frac{2\mathcal P_n}{n+1}
\left(\int_{K_0}|f_Q(\bn)|^p\,\d\omega_n(Q)\right)R^{-p}
+O_{n,f,\varepsilon}(R^{-p-2+\varepsilon}).
\end{align}
By \eqref{SLaction}, we have that $|f_{-Q}(\vecw)|=|f_Q(-\vecw)|$, so
\begin{equation}\label{eq:38}
2\int_{K_0}|f_Q(\bn)|^p\,\d\omega_n(Q)
=\int_{\Un}|f_Q(\bn)|^p\,\d\omega_n(Q).
\end{equation}
Furthermore, the main term is independent of $\vecm\in\calE_n^0$. Using that $|\calE_n^0| = 2^n$, and $|\calE_n|=2^{n-1}(2^n+1)$, we therefore obtain \eqref{eq:generic-tail}:
  \begin{equation}
    \label{eq:37}
    \Psi_f(R) = \frac{2}{(2^n + 1)}\frac{\mathcal{P}_n}{(n + 1)\mathcal{V}_n} \left(\:\int_{\Un} |f_Q(\bn)|^{2(n + 1)}\de \omega_n(Q) \right) R^{-2(n +1)} + O_{n,f,\veps}(R^{-2(n + 2) + \veps}).
\end{equation}
Positivity of the $f$ dependent constant follows from Lemma~\ref{lem:generic-positivity}.
\end{proof}

\begin{lemma}\label{lem:generic-positivity}
For $f\in\Si(\R^n)$, $\mathcal M_{n,f}>0$ if and only if $ P_+f\neq0$.
\end{lemma}

\begin{proof}
  We note that $P_+f$ is the projection of $f$ onto the space $L^2_{\text{even}}(\R^n)$, the space of even square integrable functions, which is an irreducible subrepresentation for $R$ by \cite[Theorem 4.56]{Folland-Harmonic-Analysis}. If $P_- f := f - P_+f$, we have that $P_-f$ is odd, and therefore, $R(k(Q))(P_- f)$ is odd. Therefore, $[R(k(Q))(P_- f)](\bn) = 0$, so if $P_+ f = 0$, then $\mathcal{M}_{n,f} = 0$.

  Conversely, suppose $\mathcal{M}_{n,f}=0$. Since $Q\mapsto|f_Q(\mathbf{0})|$ is continuous, $f_Q(\mathbf{0})=0$ for every $Q\in\mathrm{U}(n)$. Writing $g=n(X)a(Y)k(Q)$  gives $|(R(g)f)(\mathbf{0})|=(\det Y)^{1/4}|f_Q(\mathbf{0})|$, so $(R(g)f)(\mathbf{0})=0$ for all $g\in\mathrm{Sp}(n,\mathbb{R})$.  Now, taking $g=k(\mathrm{i}I)\,n(X)$ we obtain
  \begin{equation}
    \int_{\mathbb{R}^n}e\bigl(\tfrac12\vecw X\,\trans \vecw\bigr)
    f(\vecw)\,\mathrm{d}\vecw=0
    \qquad\text{for all }X\in\mathrm{Sym}_n(\mathbb{R}).
  \end{equation}
  Let $O \in \mathrm{O}(n)$ be orthogonal and $\vecxi \in \R^n$ such that $X = O\: \diag(\vecxi) \trans O$. Then, by a change of variables we have
  \begin{equation}
    \int_{\mathbb{R}^n}e\bigl(\tfrac12\vecw X\,\trans \vecw\bigr)
    f(\vecw)\,\mathrm{d}\vecw =
    \int_{\mathbb{R}^n}e\bigl(\tfrac12 \vecxi \cdot (w_1^2, \dots, w_n^2)\bigr)
    f(\vecw \trans O)\,\mathrm{d}\vecw.
  \end{equation}
  Let $s: (w_1, \dots, w_n) \mapsto (w_1^2, \dots, w_n^2)$, and $\de \lambda_f = f(\vecw \trans O) \de \vecw$. Then
  \begin{equation}
    \int_{\mathbb{R}^n}e\bigl(\tfrac12 \vecxi \cdot (w_1^2, \dots, w_n^2)\bigr)
    f(\vecw \trans O)\,\mathrm{d}\vecw =
    \int_{\mathbb{R}^n}e\bigl(\tfrac12 \vecxi \cdot \vecv\bigr)\de[s_*\lambda_f](\vecv),
  \end{equation}
  the Fourier transform of the push-forward measure $s_*\lambda_f$. Therefore $s_*\lambda_f = 0$. Since $f = P_+ f + P_-f$, and $P_-f$ is odd, we in fact have that
  \begin{equation}
    \int_{\mathbb{R}^n}e\bigl(\tfrac12\vecw X\,\trans \vecw\bigr)
    f(\vecw)\,\mathrm{d}\vecw = \int_{\mathbb{R}^n}e\bigl(\tfrac12\vecw X\,\trans \vecw\bigr)
    P_+f(\vecw)\,\mathrm{d}\vecw,
  \end{equation}
  and so $s_*\lambda_{P_+f} = 0$, for every $O \in \mathrm O (n)$. 

  Finally, to see that $P_+f=0$, we first note that for any odd $\psi \in C_b(\R^n)$,
  \begin{equation}
    \int_{\R^n}\psi(\vecw) P_+f(\vecw) \de \vecw = 0,
  \end{equation}
  since $P_+f$ is even. Therefore the Fourier transform of $P_+f$ is given by
\begin{equation}
  \widehat{P_+f}(\vecv)
  =\int_{\R^n}P_+f(\vecw)\,e(-\vecw\trans\vecv)\,\d\vecw
  =\int_{\R^n}\cos(2\pi\,\vecw\trans\vecv)\,P_+f(\vecw)\,\d\vecw.
\end{equation}
For $t\in\R$, let $\varphi_t(\vecv):=\cos\bigl(2\pi t\sqrt{|v_1|}\bigr)$. Then $\varphi_t(s(\vecv))=\cos(2\pi t v_1)$.
Given $\vecv \in \R^n$ we may write it as $\vecv = t \vecu$, where $t \geq 0$, and $\vecu$ is a unit vector. Let $O \in \mathrm O (n)$ be such that $\trans\vecu$ is its first column. Then the 
substitution $\vecw\mapsto\vecw \trans O$ gives 
\begin{equation}
  \int_{\R^n}\cos(2\pi \,\vecw\trans\vecv)\,P_+f(\vecw)\,\d\vecw
   =\int_{\R^n}\cos(2\pi t w_1)\,P_+f(\vecw\,\trans O)\,\d\vecw
   =\int_{\R^n}\varphi_t\,\d(s_*\lambda_{P_+f}) = 0.
 \end{equation}
Fourier inversion then gives that $P_+f=0$.
\end{proof}

\section{Tail asymptotics for diagonal \texorpdfstring{$X$}{X}}
\label{section:diagonal-asymptotics}

Let $\mathbb T_K=\{\diag(q_1,\ldots,q_n):q_j\in\mathrm U(1)\}
\subset\Un$ be the diagonal unitary torus. We use $Q\in\mathbb T_K$
for its elements and $f_Q=R(k(Q))f$ as before.
Write $\d\omega_K(Q)$ for the product Haar probability measure
$\d\omega_1(q_1)\cdots\d\omega_1(q_n)$.
Define
\begin{equation}\label{eq:diagonal-angular-moment}
\mathcal M_f^\Diag:=\int_{\mathbb T_K}|f_Q(\bn)|^4
\,\d\omega_K(Q).
\end{equation}

\begin{theorem}\label{thm:diagonal-tail}
Let $n\geq2$ and $f\in\Si(\R^n)$. As $R\to\infty$,
\begin{equation}\label{eq:diagonal-tail}
\Psi_f^\Diag(R)=C_{n,f}^\Diag R^{-4}(\log R)^{n-1}
+O_{n,f}\!\left(R^{-4}(\log R)^{n-2}\log\log R\right),
\end{equation}
where
\begin{equation}\label{eq:diagonal-tail-constant}
C_{n,f}^\Diag=\frac{4^{n-1}}{(n-1)!}
\left(\frac2\pi\right)^n\mathcal M_f^\Diag.
\end{equation}
The coefficient is positive if and only if $P_\square f\neq0$.
\end{theorem}

We use the standard product domain
\begin{equation}\label{eq:diag-fd}
\calF_n^\Diag:=\{z\in\uhp_1:-\tfrac12\leq\Re z<\tfrac12,\ |z|\geq1\}^n.
\end{equation}
Write $\vecz=\vecu+\i\vecv$, and parametrise the compact factor by
$Q(\vecphi)=\diag(\e^{\i\phi_1},\ldots,\e^{\i\phi_n})$, where $\vecphi=(\phi_1,\ldots,\phi_n)$. Up to a null set, $\calF_n^\Diag\times[0,\pi)^n$ is a fundamental domain for $(\sltz)^n\bs(\sltr)^n$. 
The set of even theta characteristics when $n = 1$ is given by
\begin{equation}\label{eq:rank-one-even-characteristics}
\calE_1=\{(0,0),(0,1),(1,0)\}.
\end{equation}
The orbit of $(\bn, \bn)$ under $\sltzn$ is $\calE_n^\Diag=\calE_1^n$, and $[(\sltz)^n:\oDel]=3^n$. Let $\Lambda_\Theta^\Diag$ denote a set of representatives for $\oDel\bs(\sltz)^n$. Then
\begin{equation}\label{eq:fd-diag}
\calF_\Theta^\Diag:=\bigcup_{\gamma\in\Lambda_\Theta^\Diag}
\gamma\calF_n^\Diag
\end{equation}
is a fundamental domain  for $\oDel$ on $\uhp_1^n$.
Its hyperbolic volume is $\pi^n$. The probability measure on
$\calF_\Theta^\Diag\times[0,\pi)^n$ is therefore
\begin{equation}\label{eq:diag-haar-in-coordinates}
\d\nu_\Theta
=\pi^{-2n}\prod_{j=1}^n\frac{\d u_j\,\d v_j\,\d\phi_j}{v_j^2}
=\left(\frac2\pi\right)^n
 \prod_{j=1}^n\frac{\d u_j\,\d v_j}{v_j^2}\,
 \d\omega_K(Q).
\end{equation}
In the last expression Haar measure is restricted to the half-circle
in every factor. 
Let
$\Psi_f^\Diag(R):=\nu_\Theta\{|\Theta_f|>R\}$, for all $R\in\R$, so that $\Psi_f^\Diag(R)=1$ for $R<0$. Then
\begin{align}\label{eq:diagonal-characteristic-decomposition}
\Psi_f^\Diag(R)&=\sum_{\vecm\in\calE_n^\Diag}I_{\vecm}^\Diag(R), \\
I_{\vecm}^\Diag(R)&:=\int_{\calF_n^\Diag\times[0,\pi)^n}
 \ind_{\{|\Theta_f(\vecm/2,g)|>R\}}\,\d\nu_\Theta(g).
\end{align}

\begin{lemma}\label{lem:simplex}
For $k\geq1$ and $T\geq1$,
\begin{equation}\label{eq:Jk}
J_k(T):=\int_{\substack{v_j\geq1\ (1\leq j\leq k)\\v_1\cdots v_k>T}}
\prod_{j=1}^k\frac{\d v_j}{v_j^2}
=\frac1T\sum_{m=0}^{k-1}\frac{(\log T)^m}{m!}.
\end{equation}
\end{lemma}

\begin{proof}
We use the substitution $s_j=\log v_j$.  Then $\mathrm dv_j/v_j^{2}=\mathrm e^{-s_j}\,\mathrm ds_j$, so
\begin{equation}
  \prod_{j=1}^{k}\frac{\mathrm dv_j}{v_j^{2}}
  \;=\;\mathrm e^{-w}\,\mathrm d\boldsymbol s ,
  \qquad w:=s_1+\cdots+s_k ,
\end{equation}
and the region of integration becomes $\{\,s_j\ge 0,\ w>\log T\,\}$.  The push-forward of Lebesgue measure on $\{s_j\geq0\}$ under $w=\sum_j s_j$ has density
$w^{k-1}/(k-1)!$. Therefore
\begin{equation}
  J_k(T)\;=\;\int_{\log T}^{\infty}\frac{w^{k-1}}{(k-1)!}\,\mathrm e^{-w}\,\mathrm dw
  \;=\;\mathrm e^{-\log T}\sum_{m=0}^{k-1}\frac{(\log T)^{m}}{m!},
\end{equation}
by the standard evaluation of the upper incomplete Gamma function at an integer parameter.  This gives \eqref{eq:Jk}.
\end{proof}

Since $f$ is a Schwartz cut-off, we have that for any $(\vecx,\vecy) \in \R^n \times \R^n$, and any $(\vecz, \vecphi) \in \calF_n^\Diag \times [0,\pi)^n$

\begin{equation}\label{eq:diag-theta-upper-bound}
|\Theta_f((\vecx,\vecy);\vecz,\vecphi)|\ll_{n,f}(v_1\cdots v_n)^{1/4}.
\end{equation}
If $x_j\notin\Z$, a similar argument gives for every $A>0$, 
\begin{equation}\label{eq:shifted-diagonal-cusp-bound}
|\Theta_f((\vecx,\vecy);\vecz,\vecphi)| \ll_{n,f,A,\vecx}(v_1\cdots v_n)^{1/4}v_j^{-A}.
\end{equation}

\begin{lemma}\label{lem:diag-tail-reduction-2}
Let $n\geq2$. Uniformly for $(\vecx,\vecy)\in\R^{2n}$, $R\geq2$, $T\geq1$ and $1\leq j\leq n$, we have
\begin{align}\label{eq:diag-tail-reduction-2}
&\nu_\Theta\{(\vecz,\vecphi)\in\calF_n^\Diag\times[0,\pi)^n:
 |\Theta_f((\vecx,\vecy);\vecz,\vecphi)|>R,\ v_j<T\}\notag\\
&\qquad\ll_{n,f}R^{-4}(\log R)^{n-2} (1+\log T).
\end{align}
\end{lemma}

\begin{proof}
  By a coordinate permutation it suffices to take $j=n$. Let $b=\sqrt3/2$ and let $B_{n,f}\geq1$ denote the implied constant from \eqref{eq:diag-theta-upper-bound}. We observe that
  \begin{multline}
    \{(\vecz,\vecphi)\in\calF_n^\Diag\times[0,\pi)^n:
    |\Theta_f((\vecx,\vecy);\vecz,\vecphi)|>R,\ v_n<T\}\\ \subseteq \{v_1\cdots v_n>(R/B_{n,f})^4,\ b\leq v_n<T\}.
  \end{multline}
Lemma~\ref{lem:simplex} gives
\begin{equation}\label{eq:4B}
\int_{\substack{v_i\geq b\ (i<n)\\
v_1\cdots v_{n-1}>(R/B_{n,f})^4/v_n}}
\prod_{i=1}^{n-1}\frac{\d v_i}{v_i^2}
\ll_{n,f}\frac{v_n}{R^4}(\log R)^{n-2}.
\end{equation}
As $\Re (\vecz)$ and $\vecphi$ are constrained to compact regions, we have 
\begin{multline}
  \nu_\Theta\{(\vecz,\vecphi)\in\calF_n^\Diag\times[0,\pi)^n:
 |\Theta_f((\vecx,\vecy);\vecz,\vecphi)|>R,\ v_n<T\}\\
\ll_{n,f}R^{-4}(\log R)^{n-2}\int_b^T\frac{\d v_n}{v_n}
\ll_{n,f}R^{-4}(\log R)^{n-2}(1+\log T).  
\end{multline}
\end{proof}

\begin{lemma}\label{lem:diag-tail-reduction-1}
Let $n\geq2$. If $(\vecx,\vecy)\in\tfrac12\calE_n^\Diag$ with $\vecx\neq\bn$, then, for $R\geq2$,
\begin{equation}
\nu_\Theta\{(\vecz,\vecphi)\in\calF_n^\Diag\times[0,\pi)^n:  |\Theta_f((\vecx,\vecy);\vecz,\vecphi)|>R\} \ll_{n,f}R^{-4}(\log R)^{n-2}.
\end{equation}
\end{lemma}

\begin{proof}
 Without loss of generality, we may assume $x_n\neq0$, since coordinate permutations replace $f$ by one of finitely many Schwartz functions, so the constants remain uniform. For any fixed $A >1/4$, we then have
  \begin{equation}
|\Theta_f((\vecx, \vecy), \vecu + \i\vecv, \vecphi)|\ll_{n,f,A}(v_1\cdots v_n)^{1/4}v_n^{-A} \ll  (v_1 \cdots v_{n-1})^{1/4},
\end{equation}
uniformly in $(\vecx, \vecy) \in \tfrac12\calE_1^n$ with $x_n\neq 0$. Applying Lemma~\ref{lem:simplex} with $k=n-1$, after rescaling the lower bounds $v_j\geq\sqrt3/2$, and integrating $v_n^{-2}$ over $[\sqrt3/2,\infty)$ gives
\begin{equation}
  \nu_{\Theta}(\{ (\vecu + \i\vecv, \vecphi) \in \calF^\Diag_n\times [0,\pi)^n : |\Theta_f((\vecx,\vecy);\vecu+\i\vecv,\vecphi)| > R \}) \ll_{n,f,A} R^{-4} (\log R)^{n-2}.
  \end{equation}
\end{proof}

 Finally, we observe that if $\vecx=\bn$, and $(\vecz, \vecphi) \in \calF_n^\Diag \times [0, \pi)^n$, then 
\begin{equation}\label{eq:diagonal-cusp-expansion}
\left||\Theta_f((\bn,\vecy);\vecz,\vecphi)|
-(v_1\cdots v_n)^{1/4}|f_{Q(\vecphi)}(\bn)|\right|
\ll_{n,f,A}(v_1\cdots v_n)^{1/4}\sum_{j=1}^n v_j^{-A}.
\end{equation}
This observation allows us to reduce the proof of  Theorem \ref{thm:diagonal-tail} to an integral over $v_1, \dots, v_n$.

\begin{proof}[Proof of Theorem \ref{thm:diagonal-tail}]
Fix $A\geq1$, take $T=\log R$, and let $R$ be sufficiently large. Arguing as in the proof of Theorem \ref{thm:generic-tail}, the two preceding lemmata reduce the problem, with an error $O_{n,f}(R^{-4}(\log R)^{n-2}\log\log R)$, to the $2^n$ characteristics with $\vecx=\bn$ and the region $v_j>T$ for every $j$.

Define
\begin{equation}
a_-(Q)=\max\{|f_Q (\bn)|- C_{n,f,A}T^{-A},0\},\qquad
a_+(Q)=|f_Q (\bn)|+C_{n,f,A}T^{-A}.
\end{equation}
For $\vecx=\bn$ and $v_j>T$ for every $j$, \eqref{eq:diagonal-cusp-expansion} gives, with $Q=Q(\vecphi)$,
\begin{equation}
\{(v_1\cdots v_n)^{1/4}a_->R\}
\subseteq\{|\Theta_f|>R\}
\subseteq\{(v_1\cdots v_n)^{1/4}a_+>R\}.
\end{equation}

For a fixed $\vecphi$, we apply Lemma~\ref{lem:simplex}  with the rescaling $v_j=Tw_j$. Setting $a=a_\pm(Q)>0$, we then have
\begin{equation}
\int_{\substack{v_j>T\ (1\leq j\leq n)\\v_1\cdots v_n>(R/a)^4}}
\prod_{j=1}^n\frac{\d v_j}{v_j^2}
=\frac{a^4}{R^4}\sum_{r=0}^{n-1}
\frac{(4\log R-4\log a-n\log T)^r}{r!}.
\end{equation}
For large $R$, the rescaled threshold is at least $1$, uniformly for
bounded $a$. The integral is zero at $a=0$. We point out that since $a^4|\log a|^r$ extends continuously to $a=0$ and is bounded when $a$ is bounded, the  expansion of this polynomial is uniform in $\vecphi$ even where $a_\pm$ vanishes. Thus
\begin{align}\label{eq:diag-integral-fixed-characteristic}
&\nu_\Theta\{(\vecz,\vecphi)\in\calF_n^\Diag\times[0,\pi)^n:
(v_1\cdots v_n)^{1/4}a_\pm(Q(\vecphi))>R,\ v_j>T\ 1 \leq  j \leq n\}
\notag\\
&=\frac{4^{n-1}}{\pi^n(n-1)!}
\left( \frac{1}{\pi^n}\int_{[0,\pi)^n}a_\pm(Q(\vecphi))^4\,\d \vecphi\right)
\frac{(\log R)^{n-1}}{R^4}
+O_{n,f}\!\left(R^{-4}(\log R)^{n-2}\log\log R\right).
\end{align}
Let $E_j$ be the diagonal matrix with $-1$ in coordinate $j$ and $1$ elsewhere. Then $R(k(E_j))$ is, up to a phase, the reflection $w_j \mapsto -w_j$. This commutes with $R(k(Q))$, and fixes the origin. Therefore $|f_{QE_j}(\bn)|=|f_Q(\bn)|$. It follows that
\begin{equation}\label{eq:diagonal-angular-half-circle}
\int_{[0,\pi)^n}a_\pm(Q(\vecphi))^4\,\d\vecphi
=\pi^n\int_{\mathbb T_K}a_\pm(Q)^4\,\d\omega_K(Q).
\end{equation}

  For $R$ sufficiently large we have that
  \begin{equation}
    \int_{\mathbb T_K}|a_\pm(Q)|^4 \de \omega_K(Q)  = \int_{\mathbb T_K}|f_Q (\bn)|^4 \de \omega_K(Q) + O_{n,f,A}(T^{-A}).
  \end{equation}

  Since $A\geq1$ and $T=\log R$, the resulting error
$O_{n,f,A}(R^{-4}(\log R)^{n-1}T^{-A})$ is absorbed by the claimed
error term. Finally, we observe that the main term is independent of the choice of characteristic $(\vecx, \vecy)$ with $\vecx = \bn$. Summing over these $2^n$ choices of characteristics then yields 
  \begin{multline}
    \Psi^\Diag_f(R) = 
  \lp\frac{1}{\pi}\rp^n \frac{4^{n-1}2^n}{(n-1)!} \lp \:\int_{\mathbb T_K}|f_Q (\bn)|^4 \de \omega_K(Q)\rp \frac{(\log R)^{n-1}}{R^4} \\  + O_{n,f}\lp\frac{\log \log R}{R^4} (\log R)^{n - 2}\rp,
\end{multline}
as $R\to\infty$. Positivity follows from Lemma~\ref{lem:angular-positivity}.
\end{proof}

\begin{lemma}\label{lem:angular-positivity}
For $f\in\Si(\R^n)$, $ \mathcal M_f^\Diag>0$ if and only if $P_\square f\neq0$.
\end{lemma}

\begin{proof}
This is analogous to the proof of Lemma \ref{lem:generic-positivity}. 
\end{proof}

\end{document}